\documentclass[preprint, 3p]{elsarticle}
\journal{MSSP}
\usepackage{graphicx}           % standard LaTeX graphics tool
\usepackage{amssymb,amsmath}
\usepackage{xspace}
\usepackage{color}
\usepackage{listings}
\usepackage{colortbl}
\usepackage{subcaption}
\usepackage{tikzinclude}
\usepackage{rotating}
\usepackage{longtable}
\usepackage[percent]{overpic}
\usepackage{enumitem}
\usepackage{mathtools}
\usepackage{pifont}

\usetikzlibrary{automata,positioning}
\usetikzlibrary{shapes.multipart}
\usetikzlibrary{shapes.geometric, arrows}
\usetikzlibrary{decorations.pathreplacing}
\newcommand{\nc}{\newcommand}
\nc{\rnc}{\renewcommand}
\nc{\bs}{\boldsymbol}
\rnc{\matrix}[2]{\left[\!\!\begin{array}{#1}
	#2\end{array}\!\!\right]}
\rnc{\vector}[1]{\matrix{c}{#1}}
\nc{\mm}[1]{\boldsymbol{#1}}
\nc{\mms}[1]{\boldsymbol{#1}}
\nc{\real}[1]{\Re\lbrace #1 \rbrace}
\nc{\imag}[1]{\Im\lbrace #1 \rbrace}
\nc{\kron}{\otimes}
\nc{\dd}{\mathrm{d}}
\nc{\ii}{\mathrm{i}}
\nc{\ee}{\mathrm{e}}
\nc{\inv}{^{-1}} %not used
\nc{\herm}{^{\mathrm H}}
\nc{\tra}{^{\mathrm T}}
\nc{\conj}[1]{ \overline{#1} }
\nc{\normal}{\mathrm n}
\nc{\tangential}{\mathrm t}
\nc{\kn}{{k_{\normal}}}
\nc{\kt}{{k_{\tangential}}}
\nc{\ie}{i.\,e.\xspace}
\nc{\eg}{e.\,g.\xspace}
\nc{\vs}{vs.\xspace}
\nc{\cf}{cf.\,\xspace}
\nc{\myquote}[1]{`#1'}
\nc{\etal}{et al.\xspace}
\nc{\fabstand}{\,}
\nc{\fp}{\fabstand.}
\nc{\fk}{\fabstand,}
\nc{\tab}[5][tbh]{\begin{table}[#1]\centering\caption{#4\label{tab:#5}}\begin{tabular}{#2}\hline #3 \\ \hline\end{tabular}\end{table}}
\nc{\ns}{N_\text{s}}
\nc{\fig}[4][tbh]{
\begin{figure}[#1]
\centering
\includegraphics[width=#4\textwidth]{figures/#2}
\caption{#3\label{fig:#2}}
\end{figure}}
\nc{\e}[2]{\begin{equation} #1 \label {eq:#2} \end{equation}}
\nc{\est}[1]{\begin{equation*} #1 \end{equation*}}
\nc{\ea}[1]{
\begin{align}
#1 \end{align}}
\nc{\east}[1]{
\begin{eqnarray*}
#1 \end{eqnarray*}}
\nc{\fref}[1]{{Fig.~\ref{fig:#1}}}
\nc{\frefo}[1]{{\ref{fig:#1}}}
\nc{\frefs}[1]{{Figs.~\ref{fig:#1}}}
\nc{\tref}[1]{{Tab.~\ref{tab:#1}}}
\nc{\trefo}[1]{{\ref{tab:#1}}}
\nc{\trefs}[1]{{Tab.~\ref{tab:#1}}}
\nc{\eref}[1]{{Eq.~(\ref{eq:#1})}}
\nc{\erefo}[1]{(\ref{eq:#1})}
\nc{\erefs}[1]{{Eqs.~(\ref{eq:#1})}}
\nc{\sref}[1]{{Section~\ref{sec:#1}}}
\nc{\srefo}[1]{\ref{sec:#1}}
\nc{\srefs}[1]{{Sections~\ref{sec:#1}}}
\nc{\ssref}[1]{{Subsection~\ref{sec:#1}}}
\nc{\ssrefo}[1]{\ref{sec:#1}}
\nc{\ssrefs}[1]{{Subsections~\ref{sec:#1}}}
\nc{\aref}[1]{{{\ref{asec:#1}}}}
\nc{\inst}[1]{$^{#1}$}
\nc{\naft}{N_{\mathrm{AFT}}}
\nc{\ndof}{N_{\mathrm{DOF}}}
\nc{\Neq}{N_{\mathrm{eq}}}
\nc{\nfact}{N_{\mathrm{fact}}}
\nc{\nsolpt}{N_{\mathrm{pt}}}
\nc{\nord}{P}
\nc{\nnewtavg}{\overline{N}_{\mathrm{newt}}}
\nc{\regeps}{\varepsilon_{\mathrm{reg}}}
\nc{\epstol}{\varepsilon_{\mathrm{tol}}}
\nc{\sign}{\operatorname{sgn}}

\makeindex

\begin{document}

\begin{frontmatter}
\title{
Are Chebyshev-based stability analysis and Urabe's error bound useful features for Harmonic Balance?
}
%\title{Nonlinear dynamics of self-excited friction-damped systems close to internal resonance}
\author{
Lukas Woiwode (\ding{41})$^a$\\
Malte Krack$^b$}
\address{$^a$University of Stuttgart, Stuttgart, Germany, e-mail: Lukas.Woiwode@ila.uni-stuttgart.de}
\address{$^b$University of Stuttgart, Stuttgart, Germany}
\begin{abstract}
Harmonic Balance is one of the most popular methods for computing periodic solutions of nonlinear dynamical systems.
In this work, we address two of its major shortcomings:
First, we investigate to what extent the computational burden of stability analysis can be reduced by consistent use of Chebyshev polynomials.
Second, we address the problem of a rigorous error bound, which, to the authors' knowledge, has been ignored in all engineering applications so far.
Here, we rely on Urabe's error bound and, again, use Chebyshev polynomials for the computationally involved operations.
We use the error estimate to automatically adjust the harmonic truncation order during numerical continuation, and confront the algorithm with a state-of-the-art adaptive Harmonic Balance implementation.
Further, we rigorously prove, for the first time, the existence of some isolated periodic solutions of the forced-damped Duffing oscillator with softening characteristic.
We find that the effort for obtaining a rigorous error bound, in its present form, may be too high to be useful for many engineering problems.
Based on the results obtained for a sequence of numerical examples, we conclude that Chebyshev-based stability analysis indeed permits a substantial speedup.
Like Harmonic Balance itself, however, this method becomes inefficient when an extremely high truncation order is needed as, e.g., in the presence of (sharply regularized) discontinuities.
\end{abstract}


\end{frontmatter}

%\nc{\myempty}{\Box}
%
%\begin{center}
%\begin{tabular}[t]{p{1.1cm}p{6.2cm}p{1.0cm}p{6.3cm}}
%\hline
%\multicolumn{4}{l}{\textbf{Symbols and abbreviations}}\\
%	       \multicolumn{2}{l}{\textsl{Latin letters}} 	&	       \multicolumn{2}{l}{\textsl{Greek letters}}                  \\	
%       \multicolumn{2}{l}{\textsl{Sub-, superscripts, operators}}                                       			&	               	&	        \\
%              	&	          	&	              \multicolumn{2}{l}{\textsl{Abbreviations}}                  \\
%\hline
%\end{tabular}
%\end{center}

%%========================================================================================================================
%% Introduction
%%========================================================================================================================
\section{Introduction}
% HARMONIC BALANCE + ITS POPULARITY
Harmonic Balance (HB) is a method for the computation of periodic solutions of nonlinear ordinary and differential-algebraic equation systems.
HB is a Fourier-Galerkin method where an approximation is sought in the form of a truncated Fourier series, and the corresponding residual term of the differential equation system is required to be orthogonal to the Fourier base functions contained in the ansatz \cite{urab1965}.
The resulting algebraic equation system is commonly solved using a Newton-type method in conjunction with a numerical path continuation technique.
To treat generic nonlinearities conveniently within HB, the Alternating Frequency-Time scheme \cite{came1989} is typically used, where the Fourier coefficients of the nonlinear terms are computed using the discrete Fourier transform.
Using HB, the simulation of the possibly long transient is avoided by directly determining the periodic limit state.
Near resonance, reasonable accuracy is often achieved already with a few Fourier terms.
HB reduces the computational effort often by 2-4 orders of magnitude compared to numerical integration, especially for stiff problems{, see \cite{padm1995b} (ch. 4.2) and \cite{Krack.2019} (p. 115)}. %\cite{padm1995b,Krack.2019}.
This makes HB a popular method not only in structural dynamics, but also fluid dynamics \cite{Hall.2002} and electrical circuit analysis \cite{Gilmore.1991}.
In the following, we place the focus on nonlinear mechanical systems.
%\\
% LIMITATIONS
HB has two major shortcomings, (a) challenges with stability analysis and (b) the lack of an error bound, as detailed in the following two paragraphs.
\\
% STABILITY ANALYSIS
 {
Throughout this work, asymptotic stability in the sense of Lyapunov is considered.
}
Knowing the asymptotic stability of the computed periodic responses is crucial to distinguish physically feasible from infeasible periodic responses, and to analyze bifurcations.
Two fundamentally different approaches are known for analyzing the asymptotic stability of HB approximations, the \emph{monodromy matrix method} and \emph{Hill's method}.
For the former, the differential equation is linearized around the approximated periodic orbit, and the fundamental matrix associated with this linearized problem is computed \cite{Hsu.1972,Friedmann.1977}.
The monodromy matrix is the fundamental matrix evaluated at the end of one period.
According to the Floquet theorem, the orbit is asymptotically stable if all eigenvalues of this matrix (Floquet multipliers) lie within the unit disk in the complex plane.
%Asymptotic stability defines an asymptotic decay of small perturbations of the limit cycle, thus, all perturbed trajectories eventually lead back to the limit cycle.
%}
The main task is here to compute the monodromy matrix, which is commonly achieved by numerical time step integration.
Different integration methods can be used, and many of these have been compared by Peletan \etal \cite{peletan2013}.
They found that the Newmark Method (constant-average-acceleration variant) performed best among the selected set of methods.
In contrast to the monodromy matrix method, Hill's method is formulated entirely in the frequency domain (no time integration needed) \cite{groll2001b}.
Essentially, this yields a quadratic eigenvalue problem of dimension $d(2H+1)$ where $d$ is the number of degrees of freedom and $H$ is the harmonic truncation order.
The coefficient matrices of this eigenvalue problem are usually readily available from the solution of the HB equations.
The corresponding eigenvalues should be directly related to the Floquet multipliers.
Due to the finite truncation order in practice, however, some of the eigenvalues are more accurately approximated than others \cite{laza2010}.
Hence, there is a consensus that the eigenspectrum must be filtered for the most accurate ones, but there is no consensus on how to do that \cite{Moore.2005,laza2010}.
It is known that with all currently available filter methods, Hill's method still gives erroneous results in some cases \cite{peletan2013,Krack.2019}; hence, Hill's method is not further considered in the present work.
\\
% ERROR BOUND
As stated above, besides the problems with stability analysis, HB has another major limitation: almost nothing is known on the error of the HB approximation.
When no error bound can be given, the \emph{existence} of an exact periodic solution in the neighborhood of a HB approximation is not guaranteed.
Indeed, there are simple and frequently-studied examples, such as the Duffing oscillator with softening characteristic, where low-order HB suggests periodic response regimes that turn out not to exist \cite{Wagner.2019}.
Several authors have proposed to use some norm of the time-domain \emph{residual}, which is obtained when the HB approximation is substituted back into the initial differential equation system, as \emph{error indicator}, see \eg \cite{Wagner.2019,Wagner.2018,Grolet.2012}.
It was proposed to use such residual norms to define rules for automatically adjusting the harmonic truncation order (\emph{adaptive HB} \cite{Grolet.2012,jaum2010}).
An alternative to using the residual itself (effectively: violation of the dynamic force balance) as error indicator, is to compute and use the linear response to the non-balanced higher harmonics of the nonlinear forces
\cite{Ferri.2009,Grolet.2012}.
The underlying assumption of both error indicators is a good correlation between the time-domain residual and the \emph{error} (deviation between exact response and HB approximation).
However, there is only a poor correlation in some cases, as shown \eg for a system with dry friction nonlinearity in \cite{Ferri.2009}.
It is remarkable that a mathematically rigorous a posteriori error bound was already developed 1965 by Urabe \cite{urab1965}. {\footnote{More specifically, the error bound was developed, along with existence and convergence theorems, for non-autonomous first-order ordinary differential equation systems in \cite{urab1965}. These were later generalized by Stokes \cite{Stokes.1972}.}}
However, this was so far only applied by Urabe himself \cite{urab1965,urab1966} and by van Dooren \cite{vanDooren.1988}, and only to single-degree-of-freedom oscillators with cubic nonlinear terms.
It should be remarked that alternative a posteriori methods for rigorously obtaining an error bound have been proposed, for instance those proposed by Baker \etal \cite{Baker.2005}, and Lessard \etal \cite{Lessard.2014,Lessard.2015}.
An advantage of the approach in \cite{Baker.2005} is that it requires only continuity of the vector field, whereas continuous differentiability is required in \cite{urab1965,Lessard.2014,Lessard.2015}.
A drawback of that approach is that homology computations are explicitly required, which is expected to yield higher computation effort.
The approach in \cite{Lessard.2014,Lessard.2015} is not expected to differ substantially in terms of computation effort compared to Urabe's approach.
 {
However, in the authors' opinion, Urabe's approach has the important benefit that it can be rather easily implemented in such a way that it is applicable to quite generic nonlinearities.
Thus, only Urabe's error bound \cite{urab1965} is investigated in the present work.
}
%We will focus on Urabe's error bound \cite{urab1965} in the present work.
Further, it is noteworthy that first steps towards a priori error estimation were done in \cite{Kogelbauer.2021}.
The method is still limited to systems with linear viscous damping, and will not further be considered in this work.
\\
% AUTO related:
%Apart from classical time integration, the computation of periodic solutions of ordinary differential equations is generally not limited to HB, \ie a Galerkin method considering ansatz and weight functions based on Fourier basis functions.
%For collocation, the weight functions are represented by Dirac delta distributions, while any finite-dimensional space of ansatz functions can be considered.
%A well-known choice is to select piecewise polynomial functions with Gauss-Legendre collocation points, as in the software AUTO \cite{Doedel.1999}.
%Similar to HB,ö error bounds can be estimated a posteriori \cite{deBoor.1973, Russel.1972}.
%A priori predictions allow for an efficient implementation, by adaptively choosing the density and location of collocation points \cite{Russell.1978, Auzinger.2002}.\\
% PURPOSE OF WORK AND OUTLINE
The purpose of the present work is to address the above mentioned two major shortcomings of HB.
More specifically, we seek to analyze to what extent the properties of Chebyshev polynomials can be exploited to speed up the computation of the monodromy matrix and thus the stability analysis of nonlinear dynamical systems.
 {
In particular, product and integration properties of Chebyshev polynomials permit to transform the ordinary differential equation system governing the monodromy matrix into an algebraic equation system, and thus to avoid time step integration.
In addition, Chebyshev polynomials show superior convergence behavior when approximating continuous functions \cite{boydspectral}.
%By exploiting the important properties of Chebyshev polynomials (product and integration property), it is possible to transform a differential equation into a linear algebraic equation system for the sought Chebyshev coefficients.
%In this way, a numerically expensive time step integration can be avoided.
}
%Further, we aim to determine the engineering value of Urabe's error bound for the first time.
 {
Further, we aim to determine the engineering value of a rigorous a posteriori error estimation for the first time.
For this, we will focus on Urabe's error bound \cite{urab1965}.
In this work, we will present an approach for its computational evaluation and apply it towards engineering relevant problems.
}
In \sref{periodicsolution}, the HB method, the state-of-the-art methods for stability analysis as well as the proposed Chebyshev-based method are presented.
In \sref{errorestimation_method}, Urabe's error bound is explained and a computational scheme for its evaluation is proposed.
Subsequently, the numerical performance of the Chebyshev-based stability analysis and the error bound are assessed for representative examples in \sref{numericalResultsChebyshev} and \sref{numericalResultsError}, respectively.
This article ends with concluding remarks in \sref{conclusions}.

\section{Harmonic Balance and stability analysis\label{sec:periodicsolution}}
In the following, we define the equation of motion considered throughout this work.
Then, a brief recap of Harmonic Balance (HB) is given, including its popular implementation using the Alternating Frequency-Time scheme.
Then the problem of analyzing the asymptotic stability of the found periodic oscillation is formulated.
Finally, computational methods for stability analysis are presented, including two of the most popular methods, and the aforementioned original method based on Chebyshev polynomials.

\subsection{Problem setting}
% PROBLEM SETTING
Consider nonlinear mechanical systems described by the equation of motion
\ea{\Omega^2{\mms q}^{\prime\prime} + \Omega\mms D {\mms q}^{\prime} + \mms K\mms q + \mms f_{\mathrm{nl}}(\mms q) - \mms f_{\mathrm{ex}}(\tau) = \mms 0\fk \label{eq:eqm}}
where $\mms q\in\mathbb R^{d\times 1}$ is the vector of $d$ generalized coordinates, $\square^\prime$ denotes derivative with respect to normalized time $\tau$, which is related to the non-normalized time $t$ via $\tau = \Omega t$.
The external force is $2\pi$-periodic, $\mms f_{\mathrm{ex}}(\tau+2\pi) = \mms f_{\mathrm{ex}}(\tau)$.
The vector of nonlinear forces, $\mms f_{\mathrm{nl}}$, depends only on $\mms q$.
The coefficient matrices $\mms D, \mms K \in\mathbb R^{d\times d}$ are time-invariant.
The models considered throughout the numerical examples are all of the form in \eref{eqm}.
The following formulations can be generalized, among others, to $\mms q^\prime$-dependent and explicitly time-dependent nonlinear forces, non-trivial coefficient matrices of vector ${\mms q}^{\prime\prime}$, and time-periodic coefficient matrices for different derivative orders.
For convergence of HB, the applicability of Urabe's error bound and the Floquet theorem, sufficient differentiability of $\mms f_{\mathrm{nl}}$ and $\mms f_{\mathrm{ex}}$ with respect to $\mms q$ and $\tau$, respectively, must be required \cite{urab1965} {(lemma 2.2)}.
The numerical examples are restricted to analytic functions (infinitely differentiable).

\subsection{Harmonic Balance}
% HARMONIC BALANCE
HB approximates $\mms q(\tau)$ as Fourier series, $\mms q(\tau)\approx \mms q_H(\tau)$, truncated to order $H$,
\ea{
\mms q_H(\tau) = \sum_{k=-H}^{H} \hat{\mms q}(k)\ee^{\ii k\tau}\fk \label{eq:FS}
}
which is written here in complex-exponential form with the imaginary unit $\ii = \sqrt{-1}$.
To ensure that $\mms q_H(\tau)$ is real-valued, the Fourier coefficients are pairwise complex-conjugate, $\hat{\mms q}(-k) = \overline{\hat{\mms q}}(k)$, where $\overline{\square}$ denotes complex-conjugate.
HB requires that the Fourier coefficients of the residual obtained by substituting \eref{FS} into \eref{eqm} vanish with regard to the retained harmonics (up to order $H$).
This corresponds to a Fourier-Galerkin projection and can be expressed as algebraic equation system
\ea{
\left(-k\Omega^2 + \ii k\Omega \mm D + \mm K\right)\hat{\mms q}(k) + \hat{\mms f}_{\mathrm{nl}}(k)-\hat{\mms f}_{\mathrm{ex}}(k) = \mms 0 \quad k=-H,\ldots,H\fp \label{eq:HB}
}
As already proposed by Urabe and Reiter \cite{urab1966}, the Fourier coefficients can be computed via discrete Fourier transform, \eg
\ea{
\hat{\mms f}_{\mathrm{nl}}(k) = \frac{1}{\naft} \sum_{n=0}^{\naft-1}\mms f_{\mathrm{nl}}\left(\mm q_H\left(\tau_n\right)\right)\ee^{-\ii k \tau_n} \quad \tau_n = \frac{2\pi}{\naft}n \quad k=-H,\ldots,H\fp
\label{eq:AFT}
}
This approach is now well-known as Alternating Frequency-Time scheme \cite{came1989} and is the most popular approach for evaluating the Fourier coefficients of the nonlinear terms, both for polynomial and generic nonlinearities \cite{Krack.2019} {(ch. 2.4)}.
It is common to solve the algebraic equation system in \eref{HB} with respect to the unknown Fourier coefficients of $\mms q_H$ using a Newton-type method in conjunction with analytical gradients, and to combine this technique with predictor-corrector path continuation \cite{Krack.2019} {(ch. 4.4)}.

\subsection{Stability according to Floquet theory}
% STABILITY ANALYSIS
Once a periodic oscillation, $\mm q_{\mathrm{p}}(\tau+2\pi)=\mm q_{\mathrm{p}}(\tau)$ is known (or its approximation $\mms q_H(\tau)$), the Floquet theorem is commonly used to analyze its asymptotic stability.
The idea is to consider an infinitesimal perturbation, $\Delta \mms q$, around $\mm q_{\mathrm{p}}$ and determine whether this grows over a period.
The mapping over one period is given by the monodromy matrix, which is the fundamental matrix, $\mm\Phi(\tau)$, evaluated at the end of the period, $\mm\Phi(2\pi)$.
If all eigenvalues of the monodromy matrix (Floquet multipliers) are within the unit disk in the complex plane, $\mm q_{\mathrm{p}}$ is asymptotically stable; if there is at least one eigenvalue outside the unit disk, the periodic oscillation is unstable.
By the way the Floquet multipliers cross the unit circle, the type of bifurcation can be inferred.
By linearizing \eref{eqm} for $\mms q =\mms q_{\mathrm p} + \Delta \mms q$ around $\mms q_{\mathrm p}$, one obtains a linear ordinary differential equation for $\Delta \mms q$,
\ea{
\Omega^2\Delta\mms q^{\prime\prime} + \Omega\mm D\Delta\mms q^{\prime} + \mms K \Delta\mms q + \underbrace{\left.\frac{\partial \mms f_{\mathrm{nl}}}{\partial \mms q}\right|_{\mms q_{\mathrm p}}}_{\mms J\left(\mms q_{\mathrm p}\left(\tau\right)\right)} \Delta\mms q = \mms 0\fk \label{eq:linearization}
}
where $\mms J$ is the Jacobian. % $2\pi$-periodic
It is  {useful} to cast \eref{eqm} in state-space form,
\ea{
\mm x = \vector{\mm q\\ \mm u} \quad \mm x^\prime = \vector{\mm u \\ \mm f_{\mathrm{ex}}(\tau) - \frac{1}{\Omega^2}\left(\Omega \mm D\mm u + \mm K\mm q + \mm f_{\mathrm{nl}}(\mm q)\right)} = \mm F\left(\mm x,\tau\right)\fp \label{eq:statespace}
}
The fundamental matrix $\mm\Phi$ contains as columns the solution of \eref{linearization} for unit initial conditions.
This can be compactly written as matrix initial value problem,
\ea{
\mm\Phi^{\prime} = \underbrace{\matrix{cc}{\mm 0 & \mm I\\ -\frac{1}{\Omega^2}\left(\mm K+\mm J\left(\mms q_{\mathrm p}\left(\tau\right)\right)\right) \qquad -\frac{1}{\Omega}\mm D}}_{\mm A\left(\mm x_{\mathrm p}\left(\tau\right)\right)}\mm\Phi \quad \mm\Phi(0) = \mm I\fk \label{eq:MIVP}
}
where $\mm x_{\mathrm p}(\tau)$ is the $2\pi$-periodic state space representation of the periodic oscillation, and $\mm I$ denotes identity matrices of according dimension.

\subsection{Computational stability analysis\label{sec:compstab}}
% COMPUTATIONAL (IMPLEMENTATION OF) STABILITY ANALYSIS
Three different methods for computing the monodromy matrix $\mm\Phi(2\pi)$ are outlined in the following, the \emph{Matrix Exponential Method (MExp)}, the \emph{Newmark Method (NTP)} and the \emph{Chebyshev Method (Cheby)}.
While the first two methods are commonly applied, the last one has never been used for analyzing the stability of periodic oscillations of nonlinear mechanical systems.
Chebyshev polynomials have, however, been proposed for computing monodromy matrices of \emph{linear} mechanical systems with time-periodic coefficient matrices \cite{sinh1991, Wu.1994}.

\subsubsection{Matrix Exponential Method (MExp)}
%       MATRIX EXPONENTIAL
The Matrix Exponential Method simply approximates the matrix $\mm A$ in \eref{MIVP} as piecewise constant between two time levels.
The only time-varying part of $\mm A$ is $\mm J$. %is commonly available on the grid points of the discrete Fourier transform in \eref{AFT}.
An expression for $\mm J$ is commonly available, as this is an auxiliary variable when evaluating analytical gradients (typically used in conjunction with Newton-type methods for solving the HB equations).
%This is because $\mm J$ is an auxiliary variable when evaluating analytical gradients (commonly used in conjunction with Newton-type methods for solving the HB equations).
For $\mm A(\tau) \approx \mm A(\tau_n)$ within $\tau_n\leq \tau \leq \tau_{n+1}$, it holds that $\mm\Phi(\tau_{n+1}) = \ee^{\mm A(\tau_n)~\left(\tau_{n+1}-\tau_n\right)}\mm\Phi(\tau_n)$, where $\ee^{\square}$ denotes the matrix exponential.
 {
Here and in the following, we write $\mm A(\tau)$ and $\mm J(\tau)$ short hand for $\mm A\left(\mm x_{\mathrm p}\left(\tau\right)\right)$ and $\mm J\left(\mms q_{\mathrm p}\left(\tau\right)\right)$, respectively.
}
Applying this on a regular grid,
\ea{
\tau_1=0\fk \tau_{N+1}=2\pi\fk \tau_{n+1}-\tau_n = 2\pi/N\fk \label{eq:equidistantgrid}
 }
yields
\ea{
\mm\Phi(2\pi) = \prod_{n=1}^{N} \ee^{\mm A(\tau_n)\frac{2\pi}{N}}\fp \label{eq:mexp}
}
The product in \eref{mexp} is to be carried out as successive multiplication from the left.
\\
% TEMPORAL REFINEMENT COMPARED TO HB/AFT
It seems tempting to use the grid points from the discrete Fourier transform in \eref{AFT}, and simply set $N=\naft$.
This is because $\mm J$ is readily available upon solution of the HB equations, as mentioned above.
However, it turns out that commonly a much finer temporal resolution is needed for stability analysis.
For a cubic-order nonlinearity, for instance, a number of $\naft =4H+1$ samples is sufficient to avoid aliasing within the Alternating Frequency-Time scheme \cite{Woiwode.2020} {(appx. A)}, whereas one to two orders of magnitude more time steps may be required to achieve reasonable accuracy of the Floquet multipliers, as will be shown later.
In practice, therefore, $\mm J$ has to be evaluated on a finer grid with $N>\naft$.
This applies also to the method subsequently described.

\subsubsection{Newmark Method (NTP)}
%       NEWMARK METHOD
The Newmark Method for computing the monodromy matrix directly uses \eref{linearization} as point of departure.
More specifically, the constant-average-acceleration variant is used for numerical integration.
Again, the grid defined in \eref{equidistantgrid} is used with the constant (normalized) time step $2\pi/N$.
This leads to the scheme:
\ea{
\mm S_{n+1} \mm\Delta \mm q_{n+1} &= \mm b_{n+1}\fk \label{eq:NTPone} \\
\Delta\mm u_{n+1} &= \frac N\pi\left( \Delta\mm q_{n+1}-\Delta\mm q_n\right) - \Delta\mm u_n\fk \label{eq:NTPtwo}\\
\Delta\mm u^\prime_{n+1} &= \frac{N^2}{\pi^2}\left( \Delta\mm q_{n+1}-\Delta\mm q_n \right) - \frac{2N}{\pi}\Delta \mm u_n - \Delta\mm u^\prime_n\fk \label{eq:NTPthree} \\
\text{with}~~ \mm S_{n+1} &= \left(\frac{N\Omega}{\pi}\right)^2\mm I + \frac{N\Omega}{\pi} \mm D + \mm K + \mm J\left(\tau_{n+1}\right)\fk\\
\mm b_{n+1} &=\Omega^2\left(\frac{N^2}{\pi^2}\Delta \mm q_n + \frac{2N}{\pi}\Delta\mm u_n + \Delta\mm u^\prime_n\right) + \Omega\mm D\left(\frac{N}{\pi}\Delta\mm q_n + \Delta\mm u_n\right)\fp
}
As in \eref{MIVP}, the scheme is applied in matrix form to the set of unit initial values, $\Delta \mms q_1 = \matrix{cc}{\mm I & \mm 0}$, $\Delta \mms u_1 = \matrix{cc}{\mm 0 & \mm I}$ which are each of dimension $d\times 2d$.
To evaluate \eref{NTPone} at $n=1$, $\Delta\mm u^\prime_1 = \Delta\mm q^{\prime\prime}_1$ is required, which is obtained with the initial values and $\mm J\left(\tau_1\right)$ from \eref{linearization}.
The scheme in \erefs{NTPone}-\erefo{NTPthree} is then applied for $n=1,\ldots,N$, to obtain $\mm\Phi(2\pi) = [\Delta\mm q_{N+1}; \Delta\mm u_{N+1}]$.

\subsubsection{Chebyshev Method (Cheby)}
%       CHEBYSHEV METHOD
The idea of the Chebyshev Method is to seek an approximate solution of \eref{linearization} in the form of a Chebyshev polynomial.
One should recall that although the considered oscillation, $\mm q_{\mathrm p}$, is $2\pi$-periodic, and hence \eref{MIVP} has a $2\pi$-periodic coefficient matrix, the sought matrix, $\mm\Phi(\tau)$, is, of course, not periodic.
Consequently, it is not possible to approximate $\mm\Phi(\tau)$ as a Fourier polynomial, but one has to resort to more generic polynomials such as the Chebyshev polynomials.
Chebyshev polynomials have many favorable properties, including the ability to express integration and multiplication using operational matrices, \cf \aref{Cheby}.
This permits to derive a linear algebraic equation system for the sought Chebyshev coefficients.
In this work, we use Chebyshev base functions of the first kind, $T_j(\tau)$, which are shifted from the standard interval $[-1,1]$ to the interval $\tau\in [0,2\pi]$.
Let us define $a_{C}(\tau)$ as Chebyshev polynomial,
\ea{
a_{C}(\tau) = \sum_{j=1}^{C} T_j(\tau) \check{a}(j) \quad \tau \in [0,2\pi]\fk \label{eq:chebydef}
}
with the truncation order $C$ and the Chebyshev coefficients $\check{a}(j)$.
For compact notation, the coefficients are collected in a vector $\check{\mm a} = [\check{a}(1);\ldots;\check{a}(C)]$.
\\
% REPRESENTATION OF ODE AS LINEAR EQUATION SYSTEM IN CHEBYSHEV SPACE
%In Chebyshev space, integration can be expressed as a linear algebraic operation.
In the following, we derive a linear algebraic equation system for the Chebyshev coefficients of the sought approximation.
To make use of the integration rule, one first has to integrate \eref{linearization} two times, in the interval from $0$ to $\tau$,
\ea{
\Omega^2\left(\Delta\mm q^\prime\left(\tau\right)-\Delta\mm q^\prime\left(0\right)\right) + \Omega\mm D\left(\Delta\mm q\left(\tau\right)-\Delta\mm q\left(0\right)\right) + \int_{0}^{\tau} \left(\mm K+\mm J\left(\tau^\bullet\right)\right)\Delta\mm q\left(\tau^\bullet\right)\dd\tau^\bullet = \mm 0&\fk \label{intone}\\
\begin{split}
\Omega^2\left(\Delta\mm q\left(\tau\right)-\Delta\mm q\left(0\right) - \int_{0}^{\tau}\Delta\mm q^\prime\left(0\right)\dd\tau^\circ\right) + \Omega\mm D\int_{0}^{\tau}\Delta\mm q\left(\tau^\circ\right)-\Delta\mm q\left(0\right)\dd\tau^\circ& \\
+ \int_{0}^{\tau}\int_{0}^{\tau^\circ} \left(\mm K+\mm J\left(\tau^\bullet\right)\right)\Delta\mm q\left(\tau^\bullet\right)\dd\tau^\bullet\dd\tau^\circ = \mm 0&\fp \label{eq:inttwo}
\end{split}
}
Herein, some integrals are trivial because the integrand is constant; the notation as integral has the advantage that the expression can be conveniently cast into Chebyshev space (via the integration rule).
As outlined in \aref{Cheby}, integration is carried out using the operational matrix $\mm G$ in Chebyshev space.
Further, the time-dependent matrix $\mm J(\tau)$ is represented as Chebyshev polynomial, by applying the discrete Chebyshev transform described below.
Then the product of $\mm J(\tau^\bullet)\Delta\mm q(\tau^\bullet)$ in \eref{inttwo} can be carried out in Chebyshev space using the product operational matrix also defined in \aref{Cheby}.
With this, the representation of \eref{inttwo} in Chebyshev space reads:
\ea{
\mm C \Delta\check{\mm q} &= \mm c\fk \label{eq:chebylin}\\
\mm C &= \Omega^2\mm I + \Omega\mm D\kron\mm G^\mathrm{T} + \mm K\kron \left(\mm G^\mathrm{T}\right)^2 + \left(\mm I\kron \left(\mm G^\mathrm{T}\right)^2\right)\mm P\left(\check{\mm J}\right)\fk \label{eq:chebyC} \\
\mm c &= \left(\Omega^2 \mm I + \Omega\mm D\kron\mm G^\mathrm{T}\right)\Delta\check{\mm q}\left(0\right) + \Omega^2\left(\mm I\kron\mm G^\mathrm{T}\right)\Delta\check{\mm q}^\prime\left(0\right)\fp \label{eq:chebyc}
}
\eref{chebylin} is a linear algebraic equation system for the sought Chebyshev coefficients of the perturbation $\Delta\mm q$, for given initial conditions $\Delta \mm q(0)$, $\Delta \mm q^\prime(0)$.
The Chebyshev coefficients of the (time-constant) initial conditions in \eref{chebyc}, $\Delta\check{\mm q}\left(0\right)$ and $\Delta\check{\mm q}^\prime\left(0\right)$, are trivial.
Analogous to the Matrix Exponential and the Newmark Method, \eref{chebylin} is applied in matrix form to the complete set of $2d$ unit initial values.
This way, the Chebyshev coefficients of the fundamental matrix $\mm \Phi(\tau)$ are obtained.
The square coefficient matrix $\mm C$ has dimension $Cd$.
% {For more details see \eg \cite{sinh1991}.}
\\
% DISCRETE CHEBYSHEV TRANSFORM; USE OF NON-EQUIDISTANT GRID
To evaluate \eref{chebyC}, $\mm J(\tau)$ must be approximated as Chebyshev polynomial, which is to be achieved using a discrete Chebyshev transform.
As with the Matrix Exponential and Newmark Methods, $\mm J$ generally has to be re-evaluated on a finer grid.
In the case of the Chebyshev Method, it turns out that an equidistant grid leads to poor convergence due to the Runge phenomenon.
To overcome this problem, it is proposed to use a non-equidistant grid,
\ea{
\tilde\tau_n = \pi\left(1-\cos\left(\frac{n-\frac12}{C}\pi\right)\right) \quad n=\frac12,2,3,\ldots,C-1,C+\frac12\fp
\label{eq:nonequidistantgrid}
}
 {
This way, the grid becomes denser towards the interval limits according to the roots of the Chebyshev polynomials.
}
Using $n=1/2$ instead of $n=1$, and $n=C+1/2$ instead of $n=C$ ensures that the interval limits $0$ and $2\pi$ are contained.
\fref{Chebyshev_interpolation} illustrates the superior convergence behavior achieved with the proposed non-equidistant grid, compared to the case of an equidistant grid.
The interpolated functions are elements of the matrix $\mm J$ obtained for a cubic nonlinearity (a) and a regularized unilateral spring (b), \cf \sref{numericalResultsChebyshev}.
For both functions (a) and (b), the equidistant grid leads to severe numerical oscillations near the interval boundaries.
When the proposed non-equidistant grid is used, in contrast, inevitable numerical oscillations occur only near the (regularized) discontinuity.
\begin{figure*}[t]
\centering
\begin{subfigure}[c]{0.485\textwidth}
\includegraphics[width=1\textwidth]{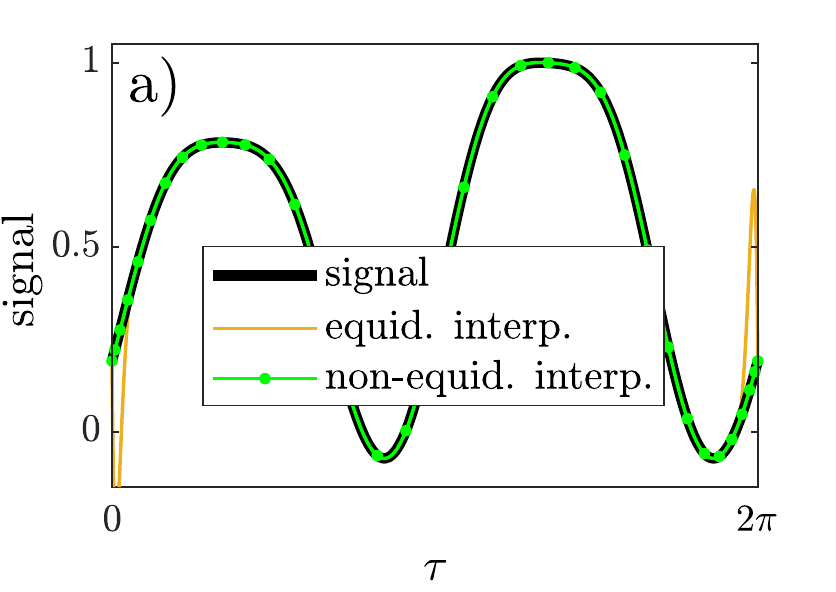}
%\subcaption{}
\end{subfigure}
\hspace{1mm}
\begin{subfigure}[c]{0.485\textwidth}
\includegraphics[width=1\textwidth]{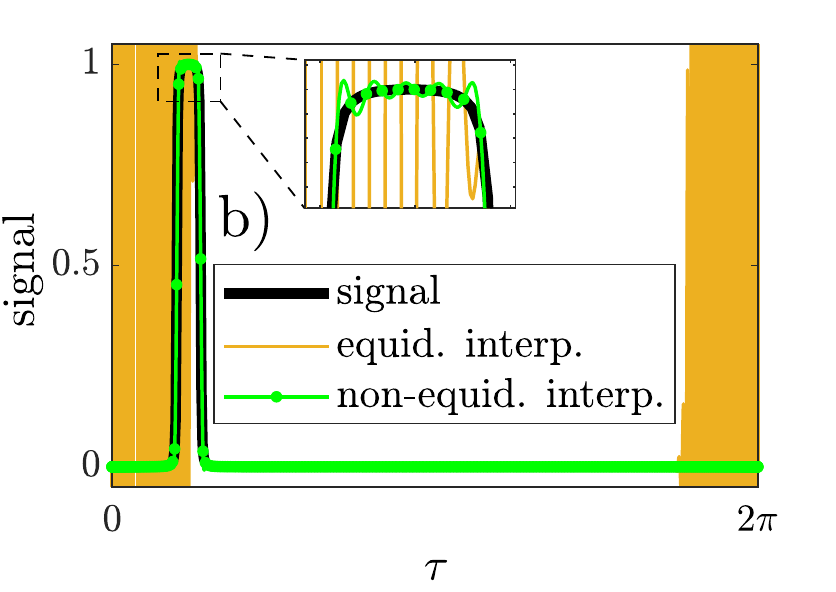}
%\subcaption{}
\end{subfigure}
\caption{Interpolation of time-periodic functions with Chebyshev polynomials using equidistant vs. non-equidistant grids: (a) function originating from the ECL benchmark (cubic nonlinearity; $C=35$); (b) function originating from oscillator with regularized unilateral spring nonlinearity ($C=300$)}
\label{fig:Chebyshev_interpolation}
\end{figure*}
\\
% IMPLICATION ON TREATMENT OF 'J'
Besides the superior convergence, the particular choice of the non-equidistant grid has the advantage that the Chebyshev polynomials are easy to evaluate.
More specifically, it is not necessary to use the recurrence definition, but one simply has $T_j(\tilde \tau_n) = \cos\left(j\left(\frac{n-\frac12}{C} - 1\right)\pi \right)$.
For the polynomial defined in \eref{chebydef}, the discrete transform can be compactly written as
\ea{
\tilde{\mm a} = \tilde{\mm T} \check{\mm a}\fk \label{eq:DTT}}
with $\tilde{\mm a} = [a(\tilde\tau_{\frac12});\ldots;a(\tilde\tau_{C+\frac12})]$ and $\mm T = [t_{nj}] = [T_j(\tilde\tau_n)]\in\mathbb R^{C\times C}$.
To obtain Chebyshev coefficients of $\mm J$, one first has to re-evaluate $\mm J$ on the non-equidistant grid, and then the inverse of $\tilde{\mm T}$ is used to solve \eref{DTT} for the Chebyshev coefficients.
Note that the re-evaluation of $\mm J$ is not an important drawback of the Chebyshev Method, as such a re-evaluation is typically needed for the other methods as well, in order to achieve reasonable accuracy.
%
%\\
%\textcolor{blue}{Check to what extent the one-time effort changes with the finding about the elements of $\mm T$ (evaluation of explicit expressions; no recurrence rule needed).}
%%\textcolor{green}{See page 7, the effort is actually slightly higher considering the explicit formulation compared to the recurrence formulae. The effort to set up $\mm T$ is, however, minor compared to setting up the other auxiliary variables/matrices.
%%So, I suggest to simply refrain from showing (also) the recursion formulae given by \eref{Chebyshev_recursive} and stick to the explicit formulation you suggested.
%%}
%\textcolor{green}{Surprisingly, applying successive vector-dot-multiplications in a single loop, as in \eref{Chebyshev_recursive}, appears to be slightly faster than determining the cosine values of a matrix spanned by a vector-vector-multiplication.}
%%
%\textcolor{green}{Therefore, we could consider to give the recursive formulation too: \\
%"(...) simply has  {$T_j(\tilde \tau_n) = \cos\left(j\left(\frac{n-\frac12}{C} - 1\right)\pi \right)$}. Yet, the recursive evaluation turns out to be numerically slightly more efficient, which can be given by
%%
%\ea{  T_{j}(\tilde \tau_n) = 2 T_{2}(\tau)  T_{j}(\tau) - T_{j-1}(\tau) , \label{eq:Chebyshev_recursive}
%}
%%
%with $T_{1}(\tilde \tau_n) = 1$ and $T_{2}(\tilde \tau_n) = \tilde \tau_n/\pi-1$."\\
%%
%However, as the one-time efforts contain the construction of multiple auxiliary variables/matrices which are much more time consuming, I suggest to refrain from showing the recurrence formulae.}

\section{Urabe's error bound and its computational implementation\label{sec:errorestimation_method}}
Urabe \cite{urab1965} established a sufficient condition for the existence of an exact periodic solution in the $\delta$-neighborhood,
\ea{
||\mms x_\mathrm{exact}- \mms x_H||<\delta\fk \label{eq:delta_neighborhood}
}
of the $H$-order HB approximation $\mm x_H$.
Here, we use the state space representation introduced in \eref{statespace}.
The condition relies on three scalar measures, an upper bound $r$ of the time-domain residual, an upper bound $\Delta$ of the variation of the Jacobian, and an upper bound $M$ for the propagation of errors within the dynamical system.
These measures are defined below and it is explained how they can be computed.
Subsequently, Urabe's existence condition is specified and the relation between $\delta$ and $r,\Delta,M$ is given.
\\
% UPPER BOUND OF TIME DOMAIN RESIDUAL
An upper bound $r$ of the time-domain residual, $\dot{\boldsymbol x}_H(\tau) - \boldsymbol F (\boldsymbol{x}_H(\tau),\tau)$, is given by:
\ea{
||\dot{\boldsymbol x}_H(\tau) - \boldsymbol F (\boldsymbol{x}_H(\tau),\tau)|| \leq || \sum\limits_{k=0}^\infty || \hat{\boldsymbol{r}}(k)||
= || \sum\limits_{k=0}^H || \hat{\boldsymbol{r}}(k)||  + \sum\limits_{k=H+1}^\infty ||\hat{\boldsymbol{r}}(k) || = r\fp
 \label{eq:timedomainresidual_estimation}
}
Herein, $ ||\cdot  ||$ denotes the Euclidian norm and $\hat{\mm r}(k)$ are the Fourier coefficients of the time-domain residual.
To determine $r$ conveniently, it is split into two parts, one associated with all harmonics up to order $H$ and one associated with higher orders.
HB requires the first part to vanish.
However, in computational practice, there is a non-zero residual, controlled by the numerical tolerance of the Newton-type solver.
This part is therefore readily available upon numerical solution of the HB equations.
As linear terms in \eref{eqm} do not generate any higher harmonics, the higher harmonic part of $r$ is only due to the nonlinear forces $\mm f_{\mathrm{nl}}$.
The higher harmonics of the nonlinear forces, generated by the HB approximation, can be computed within the Alternating Frequency-Time scheme.
For polynomial nonlinearities, one can exploit that the order of the highest harmonic is finite and can be easily determined.
For non-polynomial ones, an infinite number of higher harmonics is usually generated.
Here, a reasonable finite truncation order $H^+\gg H$ has to be selected in practice.
\\
% UPPER BOUND OF VARIATION OF JACOBIAN
$\Delta(\delta)$ is an upper bound of the variation of the matrix $\mm A = \partial \mm F/\partial \mm x$, defined in \eref{MIVP}, 
 {
\ea{
\Delta(\delta){\leq} \max_{\tau} \left|\left| ~\left. \mm A(\mm x) \right|_{\mm x_\mathrm{exact}(\tau)} ~-~ \left. \mm A(\mm x)\right|_{\mm x_{H}(\tau)}~\right|\right|,
 \label{eq:Delta_value}
}
within the $\delta$-neighborhood of $\mm x_H$ given by \eref{delta_neighborhood}.}
Since the linear forces only yield a constant part of $\mm A$, they do not contribute to $\Delta$, so only the nonlinear forces must be considered.
An analytical expression of the Jacobian $\mm J$ of the nonlinear forces is established within the Alternating Frequency-Time scheme, as explained earlier.
The variation is subsequently determined analytically for each nonlinearity individually and then accumulated to determine the Frobenius norm. 
\\
% COMPUTATION OF 'M'
 {
The third and last scalar measure required for Urabe's theorem, $M$, can be interpreted as a measure for the error propagation over one cycle. It is defined as,
}
%The third and last scalar measure, $M$, required for Urabe's theorem, is defined as
%
\begin{eqnarray}
{M}  =& \sqrt{ {2\pi} \cdot \underset{\tau}{\mathrm{max}} \left( \int\limits_0^{{2\pi}} \sum\limits_{k,l} H_{kl}^2(\tau,s)\dd s \right) }\fk \label{eq:Mvalue} %\\
\end{eqnarray}
where the matrix $\mm H = [H_{kl}]$ is given as the piecewise continuous function of the fundamental matrix $\mm\Phi$,
\ea{
\boldsymbol{H}(\tau,s)=
  \begin{cases}
    \boldsymbol{\Phi}(\tau) \left[\boldsymbol{I} - \boldsymbol{\Phi}(2\pi) \right]^{-1} \boldsymbol{\Phi}^{-1}(s)& \text{for}  \quad  0\leq s\leq \tau \leq 2\pi  \\
\boldsymbol{\Phi}(\tau)\left[\boldsymbol{I} - \boldsymbol{\Phi}(2\pi) \right]^{-1} \boldsymbol{\Phi}(2\pi)\boldsymbol{\Phi}^{-1}(s)& \text{for}\quad 0\leq \tau\leq s \leq 2\pi
\end{cases}
\fp \label{eq:rext}
}
 {
Urabe's theorem requires continuous differentiability of $\mm F(\mm x,\tau)$ and its first-order derivative in both arguments, which ensures convergence of $M$ with increasing resolution in $\tau$.}
Clearly, $M$ is computationally much more expensive to determine than $r$ or $\Delta$.
To reduce the computational burden, again, consistent use is made of the properties of Chebyshev polynomials.
From the stability analysis, the fundamental matrix $\mm\Phi(\tau$) is readily available as Chebyshev polynomial.
The multiplication and integration rules are used to first determine the matrix $\mm H$ and then to carry out multiplication and integration in \eref{Mvalue}.
Recently, it was proposed to reduce the effort for computing $M$ by using an upper bound based on the leading Floquet multipliers \cite{GarciaSaldana.2013}.
As we shall see later, a more conservative bound for $M$ will generally necessitate higher truncation orders, from which on a finite error bound can be given.
Indeed, we shall see that rather high truncation orders are necessary already with Urabe's approach, which is why the variant proposed in \cite{GarciaSaldana.2013} was not further pursued in the present work.
% It is easy to see that the computational effort required scales directly with the considered temporal resolution.
%
\begin{figure*}[h]
\centering
\includegraphics[width=0.485\textwidth]{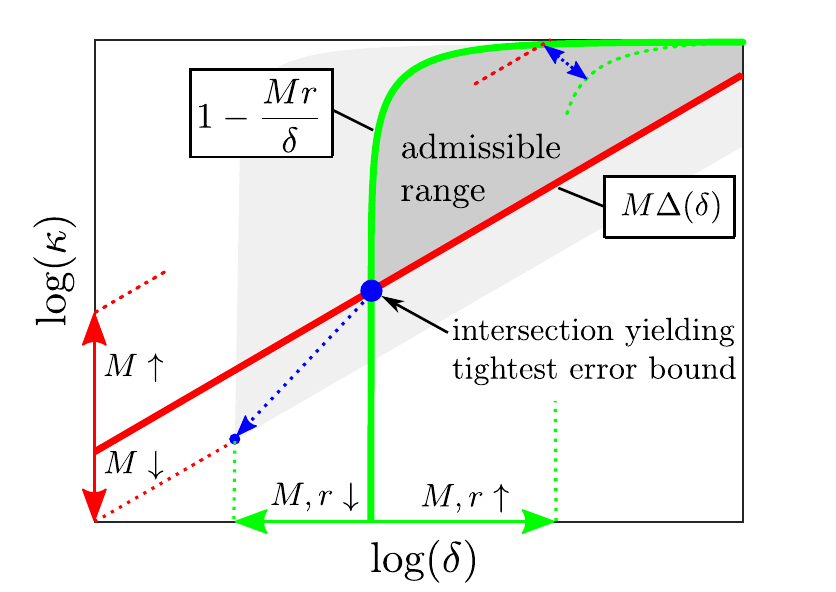}
\caption{Graphical illustration of Urabe's theorem.}
\label{fig:errorestimation}
\end{figure*}
\\
% EXISTENCE CONDITION + ERROR BOUND
Based on the three measures, $r$, $\Delta$, $M$, Urabe formulated the following theorem:
If there is a non-negative $\kappa < 1$ and a positive $\delta$ so that,
 {
\ea{\Delta (\delta)M \leq \kappa \quad \mathrm{and} \quad \kappa\leq 1 -M r / \delta,\label{eq:existencecondition}}}
then there exists one and only one periodic solution within the defined $\delta$-neighborhood.
For the considered simple examples (single-DOF oscillators with cubic nonlinearity), Urabe and Reiter \cite{urab1966} treated these inequalities analytically, in an individual way for each example.
In this work, an approach is pursued that can be more easily automatized:
From the above theorem, one can identify $M\Delta(\delta)$ as lower bound, and $1-Mr/\delta$ as upper bound for $\kappa$.
This is illustrated in \fref{errorestimation}.
One can identify the admissible range for $\kappa$ and $\delta$ in accordance with Urabe's requirement.
Obviously, the lowest $\delta$ represents the tightest error bound.
%To systematically determine this $\delta$, the intersection of upper and lower bounds for $\kappa$ is therefore calculated.
%If this intersection is within $0\leq \kappa <1$, the existence of an exact periodic solution is verified.
%The associated $\delta$ is the error bound.
In the following, we refer to the tightest error bound simply by $\delta$.
One can also see that for sufficiently large $M$ or $r$, no intersection point will occur.
In this case, no error bound for the HB approximation can be given (and hence the existence of a periodic solution is not guaranteed).
Thereby, the stronger the repelling nature of the cycle, \ie the further at least one of the Floquet multipliers lies outside the unit circle, the larger the magnitude of $M$ becomes.
%The magnitude of the $M$-value depends directly on the (in)stability of the system.
%If at least one Floquet multiplier lies far outside of the unit disk, the periodic HB solution diverges quickly to large values, which also results in large $M$-values.
%In fact, Urabe's error estimator requires that the product $Mr$ remains small to be successful.
%For large $M$, therefore,
For large $M$, $r$ must be sufficiently small so that the error estimation can be successful.
Thus, the success of the procedure can change by further increasing the truncation order $H$, as $r$ may be reduced substantially.
This way, at best, one may prove the existence of a periodic solution (and give an error bound).
However, it is impossible to prove the non-existence of a periodic solution in this way.
\\
%\COMMENT{The magnitude of the $M$-value depends directly on the (in)stability of the system.
%If at least one Floquet multiplier lies far outside of the unit disk, the periodic HB solution diverges quickly to large values, which also results in large $M$-values.
%In fact, Urabe's error estimator requires that the product $Mr$ remains small to be successful.
%For large $M$, therefore, the $r$-value must be sufficiently small so that the error estimation can be successful.}
%\COMMENT{Discuss: $M$ will depend on the stability!}

%% - INTRODUCTION:
%Besides the problems with stability analysis of HB, almost nothing is known on the error of the HB approximation.
%Urabe \cite{urab1965} developed already in 1965 a mathematically rigorous a posteriori error estimator, which hardly found attention.
%While he departed from nonautonomous ordinary differential equation systems, Stokes \cite{Stokes.1972} extended the derivations later towards autonomous systems.
%In addition to that, it proceeds from isolated periodic solutions, which applies to systems where neither of the Floquet multipliers equals one.
%Under these restrictions, it can be proven, that an HB approximation exists for sufficiently high harmonic truncation order which converges to the exact solution.
%It is further possible to establish a rigorous upper error bound to the exact solution.

\section{Numerical results: Chebyshev-based stability analysis\label{sec:numericalResultsChebyshev}}
Two numerical examples are considered in this section, the ECL benchmark involving an (inherently smooth) cubic-order stiffness nonlinearity, and a two-degree-of-freedom oscillator with a regularized elastic stop.
To this end, the methods described in \ssref{compstab} and \sref{errorestimation_method} have been implemented in MATLAB, using the open-source MATLAB tool NLvib \cite{Krack.2019} {(appx. C)} as point of departure.

\subsection{ECL benchmark\label{sec:ECL_benchmark}}
\begin{figure*}[h]
\centering
\begin{subfigure}[c]{0.325\textwidth}
\includegraphics[width=1\textwidth]{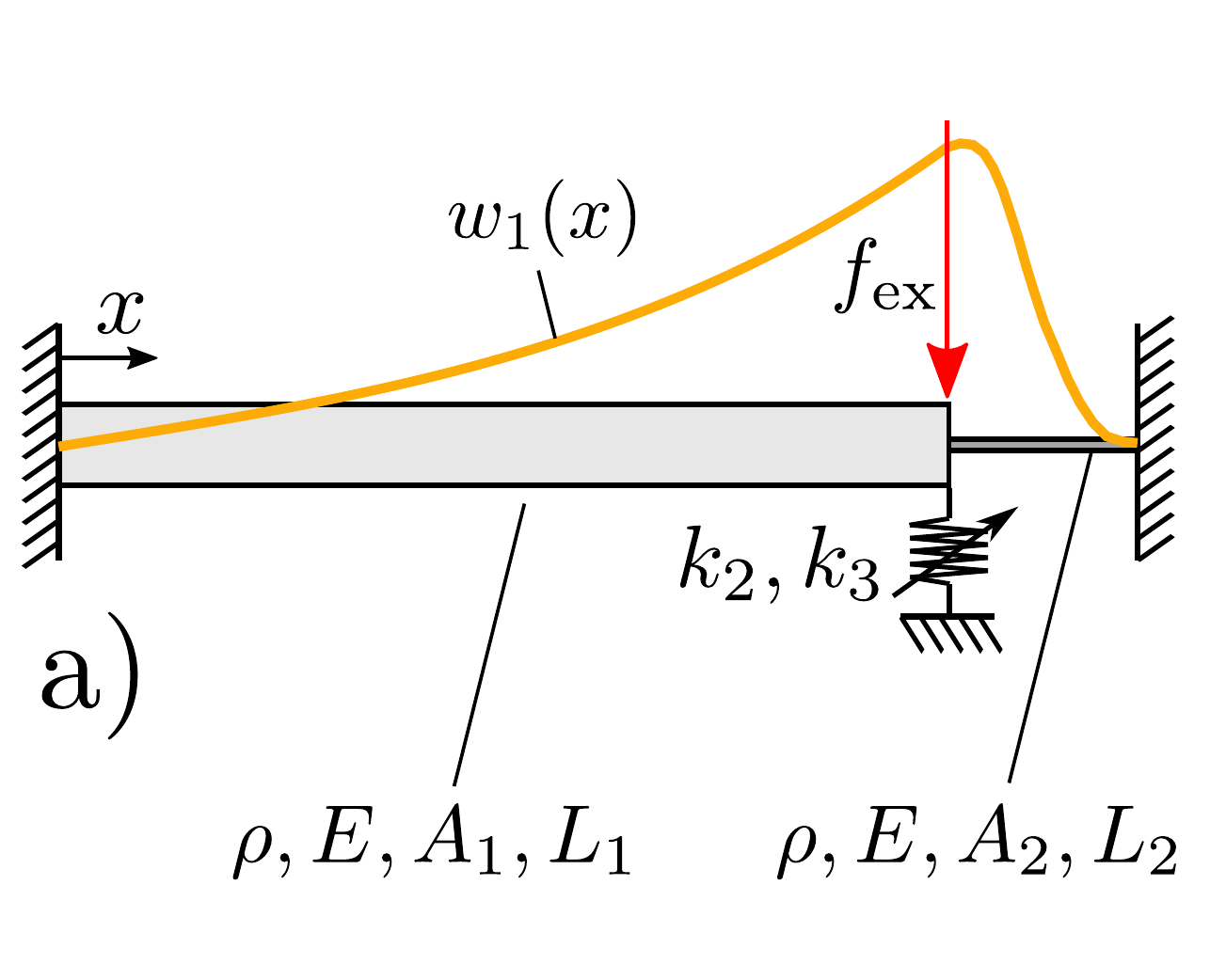}
%\subcaption{}
\end{subfigure}
\begin{subfigure}[c]{0.325\textwidth}
\includegraphics[width=1\textwidth]{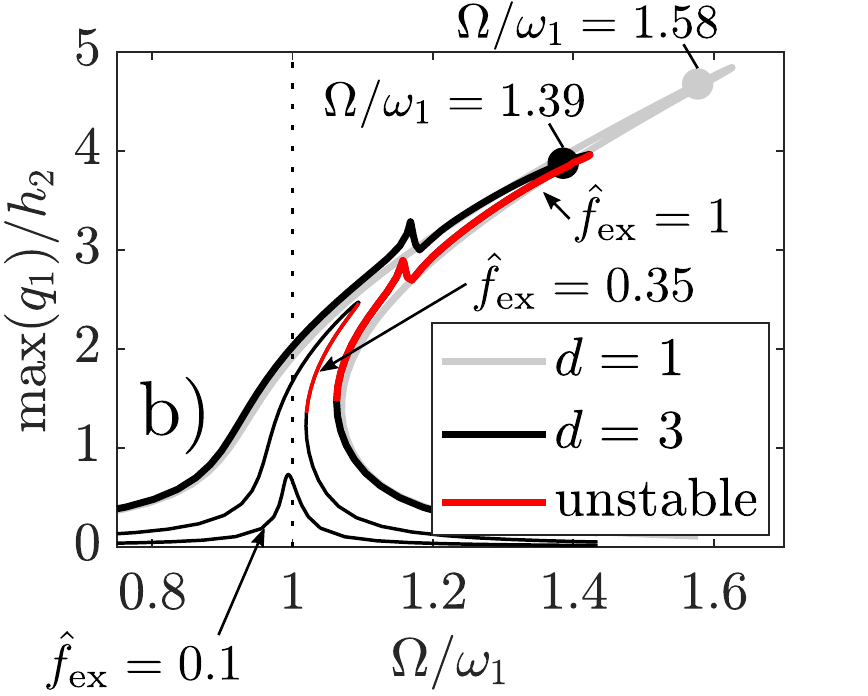}
%\subcaption{}
\end{subfigure}
\begin{subfigure}[c]{0.325\textwidth}
\includegraphics[width=1\textwidth]{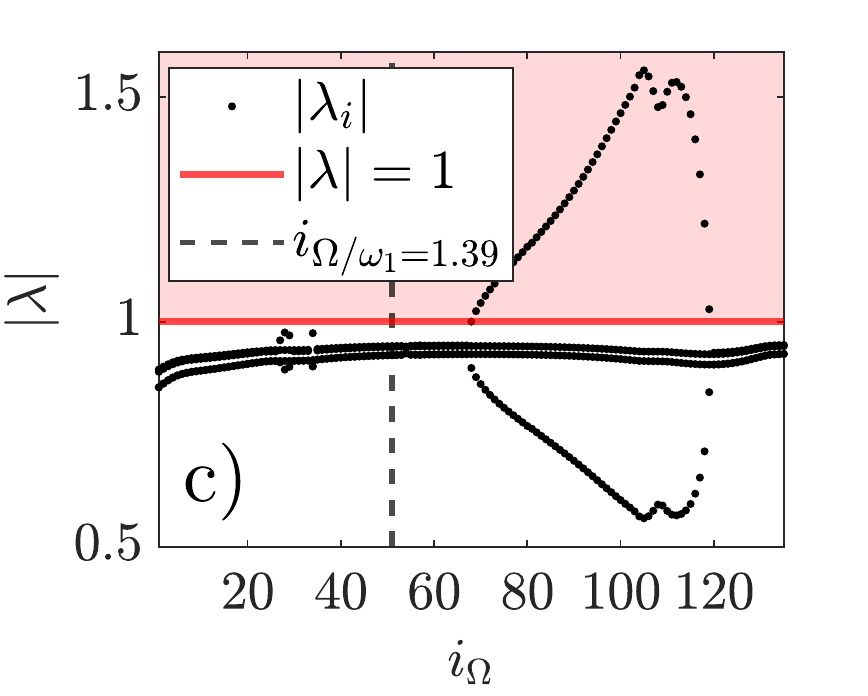}
%\subcaption{}
\end{subfigure}
\caption{ECL benchmark: (a) schematic illustration of the two connected cantilevered beams, (b) amplitude-frequency curves for different excitation levels, (c) magnitude of Floquet multipliers for the highest excitation level; material properties and geometric dimensions are density $\rho = 7800$ kg/m$^3$, Young's modulus $E=205\cdot 10^{9}$ N/m$^2$, beam cross sections $A_1=h_1 \cdot z=0.014\cdot0.014$ m$^2$ and $A_2=h_2 \cdot z=0.0005\cdot0.014$ m$^2$, beam lengths $L_1 = 0.7$ m and $L_2 = 0.04$ m
} % H = 150; Naft = 2^14; C=200; (c) is for d=3
\label{fig:ECL}
\end{figure*}
%
% DESCRIPTION OF PROBLEM SETTING
The ECL benchmark \cite{Thouverez.2003} consists of a primary beam connected to a much more slender and shorter beam, as illustrated in \fref{ECL}a.
The model and the selected parameters largely follow \cite{Noel.2015}.
The dimensions and material properties are specified in the caption of \fref{ECL}.
At the connection between the beams, a cubic spring with stiffness coefficient $k_3 = 8 \cdot 10^9 \,$N/m$^3$ and a quadratic spring with stiffness coefficient $k_2 = -1.05 \cdot 10^7 \,$N/m$^2$ is considered, modeling the hardening due to bending-stretching and the softening due to non-ideal clamping of the thinner beam, respectively.
In contrast to \cite{Noel.2015}, no finite rotational stiffness between the two beams is introduced, but instead a rigid rotational connection is modeled.
An FE model involving 1000 beam elements was set up, followed by a modal truncation to the lowest-frequency bending modes of the linearized system ($d=1$ or $d=3$).
 {
For consistency with \cite{Noel.2015}, Euler-Bernoulli elements are used, although Timoshenko elements should generally be preferred to account for shear within short beam elements.
}
Moreover, Rayleigh damping is considered, where the mass and stiffness coefficients are given by $\alpha = 3\cdot 10^{-7}$ and $\beta = 5$, respectively, resulting in a modal damping ratio of $D_1=1.2\%$ for the fundamental bending mode.
% The linear natural frequencies related to the first three bending modes are given by $\mms \omega = [209.1 , \, 944.8, \, 2610.2]$ rad/s.
\\
% RESPONSE BEHAVIOR INCLUDING STABILITY
A harmonic concentrated load, $f_\mathrm{ex} = \hat f_\mathrm{ex} \cos(\Omega t)$, is applied at the connection point.
The amplitude-frequency curves of the periodic steady-state response are depicted for different excitation levels $\hat f_{\mathrm{ex}}$ in \fref{ECL}b in the frequency range around the fundamental bending mode (natural angular frequency $\omega_1$).
The fundamental bending mode shape is also illustrated in \fref{ECL}a.
As amplitude measure, the maximum (along the period) bending displacement at the connection point is used.
A slight softening and a pronounced hardening effect can be observed, as expected for the given signs and magnitudes of $k_2$, $k_3$.
Also as expected for light damping, the amplitude-frequency curves exhibit turning points for sufficiently large excitation level.
In addition, a secondary resonance phenomenon is encountered for $\hat f_{\mathrm{ex}}=1$ in the configuration with three retained modes, $d=3$, as opposed to configuration with only one mode, $d=1$.
As expected, the overhanging branch connecting the turning points is unstable.
This is confirmed by the Floquet multipliers, whose magnitude is depicted in \fref{ECL}c along the branch obtained for $\hat f_{\mathrm{ex}}$.
Herein, $i_\Omega$ denotes the index of the solution point along the branch.
A single Floquet multiplier leaves the unit disk, as expected for turning point bifurcations.
\begin{figure*}[t]
\centering
\begin{subfigure}[c]{0.485\textwidth}
\includegraphics[width=1\textwidth]{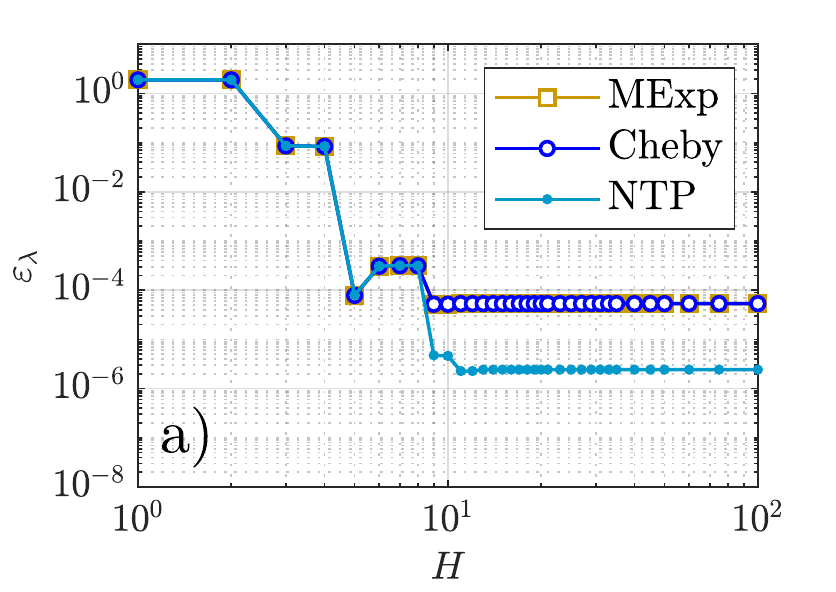}
%\subcaption{}
\end{subfigure}
\hspace{1mm}
\begin{subfigure}[c]{0.485\textwidth}
\includegraphics[width=1\textwidth]{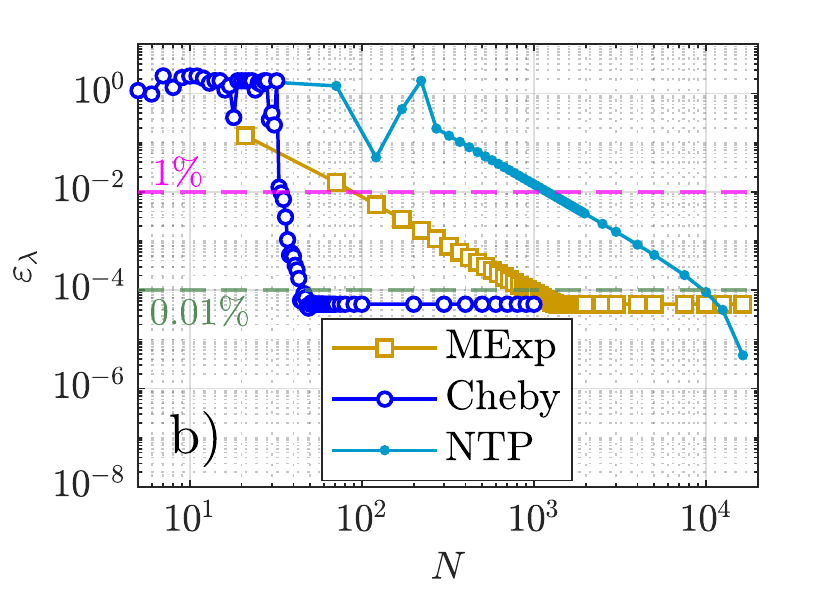}
%\subcaption{}
\end{subfigure}
%\\
%\begin{subfigure}[c]{0.485\textwidth}
%\includegraphics[width=1\textwidth]{figures/NDOF_ECL/ECL_H_Nmod3_Overview_corrected_FINERES_13-Apr-2022.pdf}
%%\subcaption{}
%\end{subfigure}
%\hspace{1mm}
%\begin{subfigure}[c]{0.485\textwidth}
%\includegraphics[width=1\textwidth]{figures/NDOF_ECL/ECL_Ntd_H9_Nmod3_Overview_corrected_FINERES_25-Apr-2022}
%%\subcaption{}
%\end{subfigure}
\caption{ECL benchmark: Leading Floquet multiplier error (a) vs. $H$, (b) vs. $N$ for the point with $\Omega/\omega_1=1.39$ ($d=3$) indicated in \fref{ECL}b
}
\label{fig:ECLstabilityHNconvergence}
\end{figure*}
\\
% H-CONVERGENCE
To assess the performance of the Chebyshev-based method of stability analysis (\emph{Cheby}), the error with respect to the leading Floquet multiplier $\lambda_{\mathrm{l}}$ is considered,
\ea{
\varepsilon_\lambda = \left|\lambda_{\mathrm{l}}-\lambda_{\mathrm{l,ref}}\right|\fk
}
where the reference $\lambda_{\mathrm{l,ref}}$ is obtained by the Shooting method in conjunction Newmark's time step integration scheme (constant-average-acceleration variant) with $N=2^{14}$ time levels per period.
The monodromy matrix $\mm\Phi(2\pi)$ is obtained as a by-product of Shooting.
The results are confronted with the Matrix Exponential (\emph{MExp}) and the Newmark Method (\emph{NTP}) of stability analysis, as described in \ssref{compstab}.
Many different solution points and model parameters (such as excitation levels and modal truncation orders) were analyzed.
Results are only shown for a representative point, namely that with $\Omega/\omega_1=1.39$ indicated in \fref{ECL}b.
The $H$-convergence is illustrated in \fref{ECLstabilityHNconvergence}a.
The number of samples within the Alternating Frequency-Time scheme was set sufficiently large to avoid aliasing errors.
The time resolution during the integration of $\mm\Phi$ was set sufficiently large to ensure that this does not have any visible effect on the depicted results.
Apparently, $H=5$ is sufficient to have an error $\varepsilon_\lambda$ smaller than $1\%$.
The individual round-off plateaus are reached at $H=9$.
The round-off plateau for the Newmark Method is slightly lower, which is attributed to the fact that essentially the same code is used as in the case of the reference method.
\\
% N-CONVERGENCE
In \fref{ECLstabilityHNconvergence}b, the $N$-convergence is illustrated for $H=9$.
Recall that generally $N\neq\naft$; in fact, $\naft=4H+1$ is sufficient to avoid aliasing for the given nonlinearity.
Thus, $N$ does not affect the HB approximation but only determines the time resolution within the stability analysis.
For the Chebyshev Method, it holds that $C=N$, as explained in \ssref{compstab}.
Apparently, the Chebyshev Method achieves high accuracy already with a coarser time resolution than the other methods.
In practice, it may be a delicate issue to properly select $C$.
It is well-known that for the task of approximating a harmonic function of order $H^*$, the required $C$ approaches $\pi H^*$ for large $H^*$,  {see \citep{boydspectral} (ch. 2.14)}.
Thus, if $H$ is the highest (non-negligible) harmonic, the cubic nonlinearity generates as highest harmonic $3H$, $C$ should be selected larger than $3\pi H$ (\eg $C=85$ for $H=9$).
Apparently, this is a good yet somewhat conservative empirical rule for the considered example.
\begin{table}[h]
\centering
%\begin{tabular}{c r l r l r l }
\caption{
 {
Computation effort of stability analysis for the error threshold $\varepsilon_\lambda<1\%$ and $\varepsilon_\lambda<0.01\%$ (offline computation effort for Cheby given in brackets): ECL benchmark with $d=1$ for $\Omega/\omega_1 = 1.58$ ($\varepsilon_\lambda<1\%$: $C=19$, $N_{\mathrm{NTP}}=N_{\mathrm{MExp}}=121$); ECL benchmark with $d=3$ for $\Omega/\omega_1 = 1.39$ ($\varepsilon_\lambda<1\%$: $C=35$, $N_\mathrm{NTP}=1221$, $N_\mathrm{MExp}=121$; $\varepsilon_\lambda<0.01\%$: $C=45$, $N_\mathrm{NTP}=10001$, $N_\mathrm{MExp}=971$); two-degree-of-freedom oscillator with elastic stop ($\varepsilon_\lambda<1\%$: $C=600$, $N_\mathrm{MExp}=1001$, $N_\mathrm{NTP}=1501$); the given computation times represent mean values based on $10,000$ repeats; solution points are indicated in \fref{ECL}b, \fref{2DOFunilateral}a, respectively. \label{tab:stability_TI_comparison} 
}
}
\begin{tabular}{c c c c c  }
 %\hline
  & {ECL benchmark ($d=1$)} & \multicolumn{2}{c}{ECL benchmark ($d=3$)} & {2DOF oscillator with elastic stop} \\
    & {$\varepsilon_\lambda< 1\%$} & {$\varepsilon_\lambda< 1\%$} & {$\varepsilon_\lambda< 0.01\%$} & {$\varepsilon_\lambda< 1\%$} \\
 \hline
 MExp   &  0.00364s	&  0.00364s   	& 0.02293s	&  0.04289s  \\
 NTP	&  0.00087s	&  0.00281s     	& 0.02593s	& 0.02700s  \\
 Cheby	&  0.00022s  &  0.00050s   	& 0.00052s 	&  0.12415s   \\
(offline)	& (0.00014s) &  (0.00034s) 	& (0.00036s)& (0.07682s)  \\
 \hline
\end{tabular}
\end{table}
\\
% COMPUTATION EFFORT
Now, we analyze the computational effort.
It is useful to recall that all methods of stability analysis need to re-evaluate the Jacobian $\mm J$ with $N<<\naft$.
Subsequently, their algorithms differ.
The Matrix Exponential Method computes a product of $N$ matrix exponentials of dimension $2d$.
The Newmark Method solves a sequence of $N$ linear equations of dimension $d$ (with $2d$-dimensional right hand side).
The Chebyshev Method solves a single linear equation of dimension $dC$ (with $2d$-dimensional right hand side), where $C$ is linked to the time resolution ($C=N$).
% RESULTS FOR FOR SIMILAR ACCURACY
%In \fref{stability_TI_comparison}a-b
In \tref{stability_TI_comparison}, the computation effort needed to meet an error threshold $\varepsilon_\lambda<1\%$ is given for $d=1$ and $d=3$, respectively, at the respective solution point indicated in \fref{ECL}b.
 {
For $d=3$, the case of a finer error threshold $\varepsilon_\lambda<0.01\%$ is additionally analyzed.
}
All computations were carried out on a Windows 10 machine with a quad-core 4GHz Intel(R) Core(TM) i7-6700K CPU, 32GB RAM  {using MATLAB R2022a}.
 {
Since the Chebyshev method converges to the reference solution for much coarser time resolution, $C<N_{\mathrm{NTP}}, N_{\mathrm{MExp}}$ holds for both error thresholds.
}
%Thanks to the quicker convergence of the Chebyshev Method, $C<N_{\mathrm{NTP}}, N_{\mathrm{MExp}}$ in both cases.
With this, the Chebyshev Method speeds up the stability analysis by about one order of magnitude 
 {
(factor $4$ and $6$ for $d=1$ and $d=3$, respectively) for $\varepsilon_\lambda<1\%$ and by about two orders of magnitude (factor $50$ for $d=3$) for $\varepsilon_\lambda<0.01\%$.
}
Here, we considered only the online effort.
For completeness, the Chebyshev Method's offline effort is additionally given in brackets in \tref{stability_TI_comparison}. %\fref{stability_TI_comparison}.
This includes setting up the matrix of the discrete Chebyshev transform, $\mm T$ and its inverse, the integration operation matrix $\mm G$ (and its products with the linear system matrices), auxiliary matrices for efficient online computation of the product operation matrix $\mm P (\cdot)$, along with the vectors of initial values.
This effort scales with the order $C$ and the size of the system matrices $d$, but can be efficiently computed once and for all (for given $d$ and maximum $C$).
Thus, this effort becomes quickly negligible during a path continuation task (with many solution points).
%In comparison, the effort for one iteration of the HB equation system ($H=5$) is $0.0012$s for $d=1$ and $0.0025$s for $d=3$, while usually between 3 and 10 iterations are needed per solution point during continuation.
%Thus, the effort required to determine the periodic solution and its stability is of same order of magnitude.
In comparison, the computation time required to solve for the periodic solution remained of the same order of magnitude as the stability estimation (one Newton iteration step with respect to the HB equations ($H=9$) took 
 {
$0.9\cdot 10^{-4}$s for $d=1$ and $2.1\cdot 10^{-3}$s for $d=3$
}
; usually 3 to 10 iterations are needed per solution point during continuation).
%\\
%\textcolor{blue}{DISCUSS ONE-TIME EFFORT (I AM ASSUMING THE INVERSION OF $\mm T$ CAUSES THE MAIN EFFORT? ENSURE THAT IT IS EXPLOITED THAT THE PROBLEM CAN BE LIMITED TO A $C\times C$ MATRIX, NOT A $dC\times dC$ MATRIX!):}
%\textcolor{green}{The inverse of $\mm T$ was already determined by the $C\times C$ matrix.}
%\\
\\
\begin{figure*}[h]
\centering
\begin{subfigure}[c]{0.485\textwidth}
\includegraphics[width=1\textwidth]{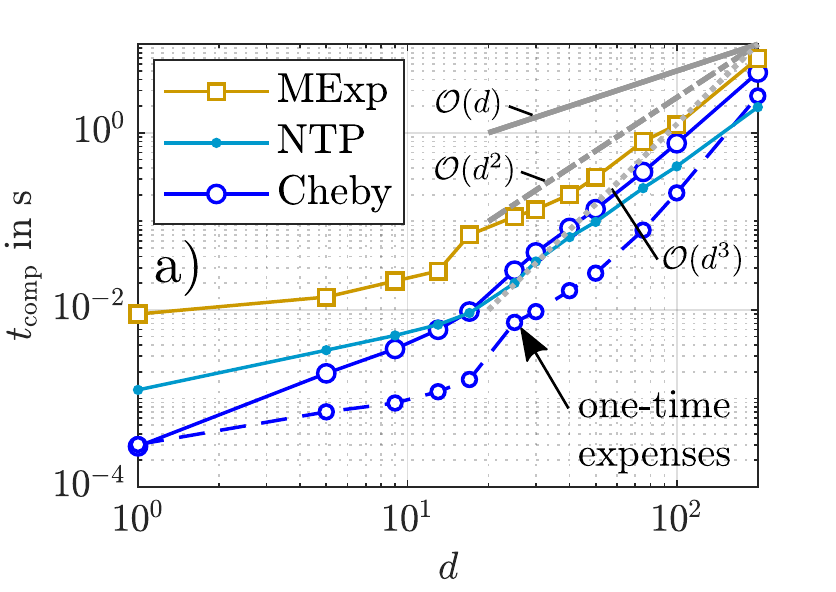}
%\subcaption{}
\end{subfigure}
\hspace{1mm}
\begin{subfigure}[c]{0.485\textwidth}
\includegraphics[width=1\textwidth]{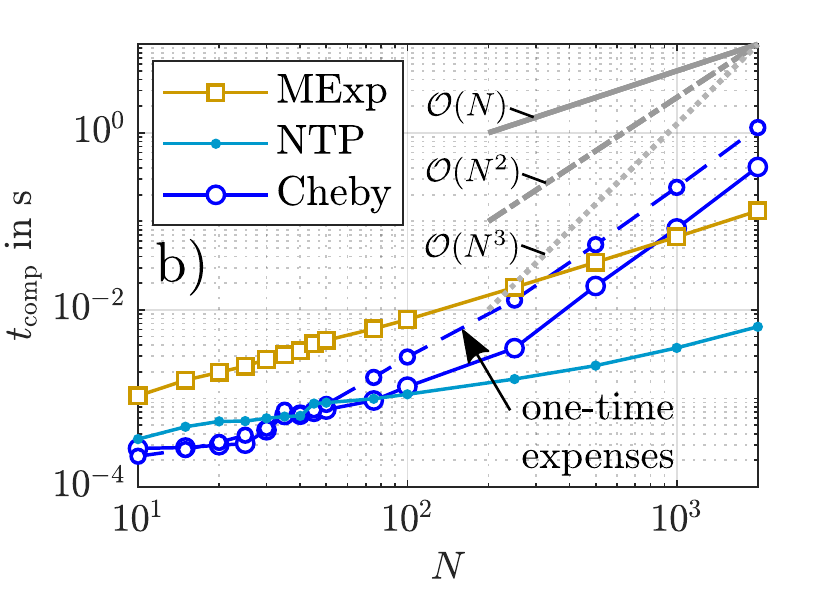}
%\subcaption{}
\end{subfigure}
\caption{Computation effort of stability analysis for ECL benchmark: (a) scaling with $d$ ($C=19$, $N_\mathrm{MExp}=N_\mathrm{NTP}=121$), (b) scaling with $N$ ($d=1$). 
 {
The depicted computation times represent mean values where the number of repeats depend on $d$ and $N$, respectively: 10,000 for the smallest values and still $d=500$ for the largest values.
}
%The computation time is normalized to the respective largest value.
}
\label{fig:stability_TI_scaling}
\end{figure*}
%
% SCALING OF COMPUTATION EFFORT WITH 'd' AND 'N'
Next, it is analyzed how the computation effort scales with the number of degrees of freedom $d$ and the time resolution $N$ (\fref{stability_TI_scaling})\footnote{
 {
The scaling behavior was computed on a machine with more memory capacity and using MATLAB R2020a, thus, the absolute values differ from values given in \tref{stability_TI_comparison}.
} 
}.
Up to about $d=17$, the Chebyshev Method shows superior computational efficiency.
Subsequently, the Newmark Method becomes more efficient, which apparently scales best among the three methods, while the Matrix Exponential Method is slowest.
Moreover, the Chebyshev Method scales worst with the time resolution.
For the given example, the Newmark Method thus becomes more efficient beyond $N=35$.
The relatively poor scaling behavior of the Chebyshev Method seems plausible in the light of the above described main computation task (solving a $dN$-dimensional equation system rather than $N$ sequential $d$-dimensional ones as in the Newmark Method).
It can thus be said that the Chebyshev Method can only be attractive, if it overcompensates its poor scaling behavior by superior convergence.

\subsection{Two-degree-of-freedom oscillator with elastic stop \label{sec:2DOFunilateral_benchmark}}
% DESCRIPTION OF PROBLEM SETTING, RESPONSE BEHAVIOR INCLUDING STABILITY
We now consider a chain of two spring-mass oscillators described by the equations of motion,
\ea{
\ddot q_1 + 0.03 \dot q_1 - 0.03\dot q_2  + q_1 - q_2 +f_\mathrm{nl}(q_1) &= 0\fk \\
\ddot q_2 - 0.03 \dot q_1 + 0.06 \dot q_2 - q_1 + 2 q_2 -0.1\cos(\Omega t) &= 0\fp
\label{eq:rext}
}
The unit mass described with coordinate $q_2$ is connected via unit springs to the ground and to the unit mass described with coordinate $q_1$.
Dashpots with damping coefficient $0.03$ are in parallel to these springs.
A harmonic forcing is applied to the second mass, while the first mass is subjected to an elastic stop (unilateral spring) with unit clearance and stiffness 100.
The nonlinear force, $f_\mathrm{nl}(q_1)$, would thus be given by $100 \cdot \mathrm{max}(q_1 - 1,0)$.
However, to ensure applicability of Urabe's and Floquet's theorems, we consider the smooth (infinitely differentiable) regularization
\ea{
f_\mathrm{nl}(q_1) = \frac{100(q_1 - 1)}{2}+\sqrt{\left(\frac{100(q_1 - 1)}{2}\right)^2 +\varepsilon_\mathrm{reg} }\fp
\label{eq:regularization}
}
The steeper the regularization, \ie the smaller $\varepsilon_{\mathrm{reg}}$, the more will \eref{regularization} resemble an elastic stop.
The effect of this regularization on the periodic response near the fundamental natural frequency is illustrated in \fref{2DOFunilateral}a.
In the following, we use $\varepsilon_{\mathrm{reg}} = 0.2$.
Besides turning points, a Torus bifurcation occurs near $\Omega/\omega_1=1$ (close to the transition from quasi-linear to nonlinear regime), and a period doubling bifurcation occurs near the resonance peak.
The magnitudes of the Floquet multipliers are depicted in \fref{2DOFunilateral}b.
\begin{figure*}[h!]
\centering
\begin{subfigure}[c]{0.475\textwidth}
\includegraphics[width=1\textwidth]{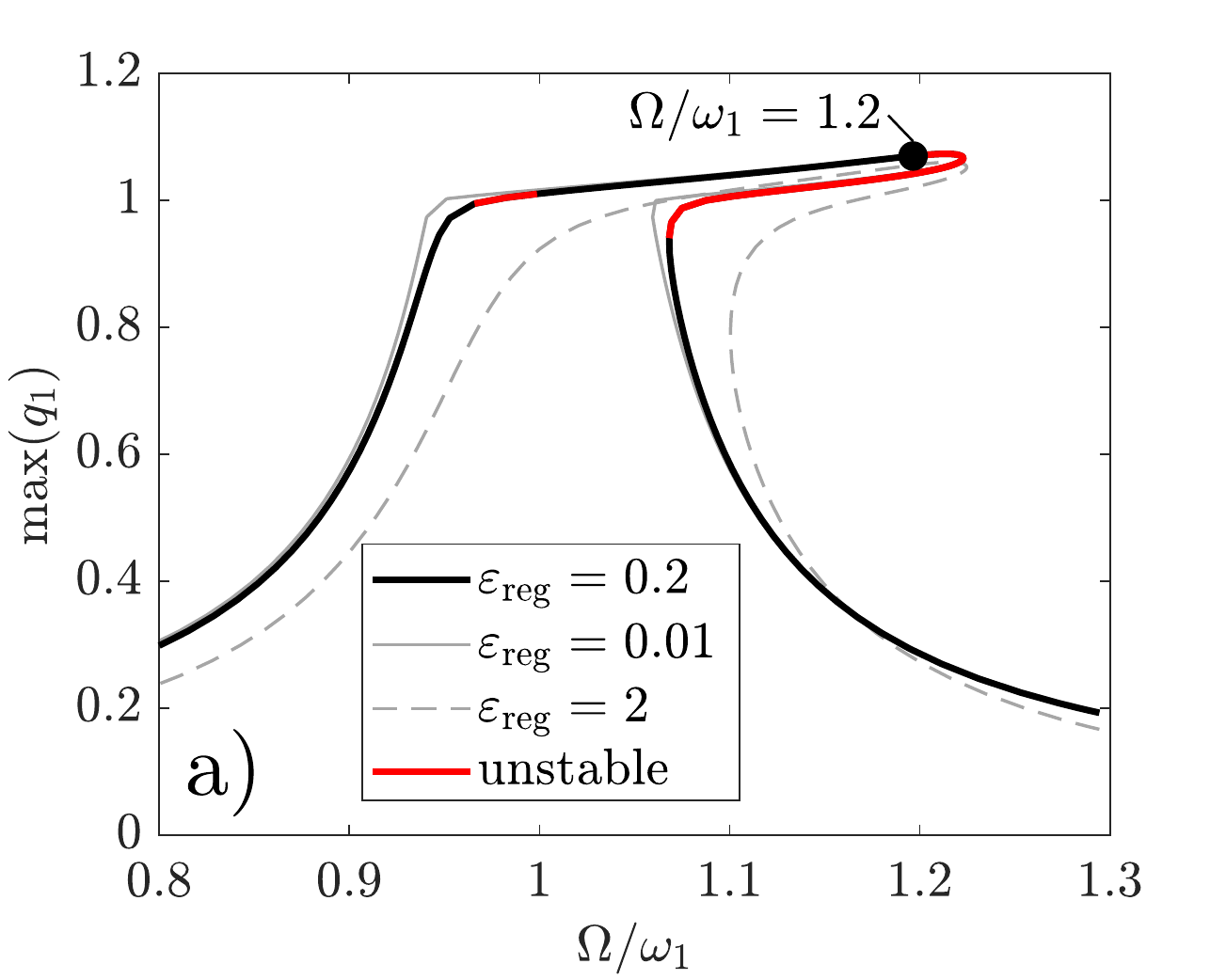}
%\subcaption{}
\end{subfigure}
\hspace{3mm}
\begin{subfigure}[c]{0.475\textwidth}
\includegraphics[width=1\textwidth]{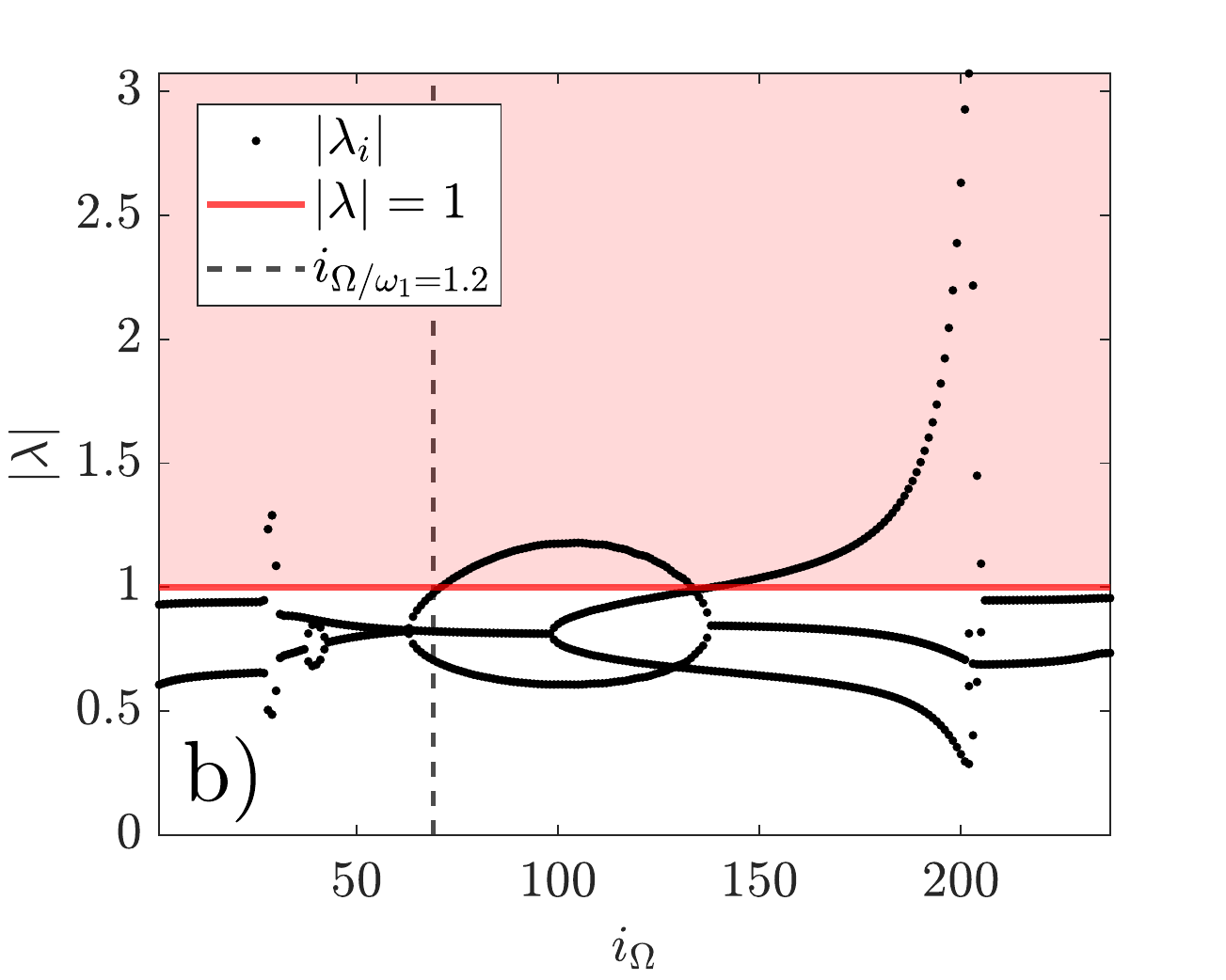}
%\subcaption{}
\end{subfigure}
\caption{Two-degree-of-freedom oscillator with elastic stop: (a) amplitude-frequency curves for different steepness of regularization $\varepsilon_{\mathrm{reg}}$, (b) magnitude of Floquet multipliers for $\varepsilon_{\mathrm{reg}}=0.2$
}
\label{fig:2DOFunilateral}
\end{figure*}
\\
% H-CONVERGENCE + COMPUTATION EFFORT
As in the previous example, many parameter sets and solution points were studied, but only representative results are shown for brevity, namely for a solution point on the upper branch near the period doubling bifurcation (indicated in \fref{2DOFunilateral}, $\Omega/\omega_1=1.2$).
The results of the $H$- and $N$-convergence are shown in \fref{2DOFunilateral_Hconv}a and b, respectively.
The reference for the definition of the leading Floquet multiplier error was again obtained using Shooting, however, a refinement to $N=2^{15}$ time levels was found to be necessary.
A much larger truncation order, {$H>200$}, is needed to reach the round-off plateau, as compared with the case of the cubic nonlinearity.
Similarly, a much finer time resolution is needed:
To achieve the specified error threshold $\varepsilon_\lambda<1\%$ {($H=80$)}, the minimum number of samples is $C=600$, $N_{\mathrm{MExp}}=1001$, $N_{\mathrm{NTP}}=1501$.
Thus, the advantage of the Chebyshev interpolation is much less pronounced for the given type of nonlinearity.
In conjunction with the poor scaling, the Chebyshev Method actually becomes inferior (see \tref{stability_TI_comparison}).
More specifically, the Chebyshev Method is  {by a factor of about $5$} slower than the Newmark Method.
\begin{figure*}[h]
\centering
\begin{subfigure}[c]{0.485\textwidth}
\includegraphics[width=1\textwidth]{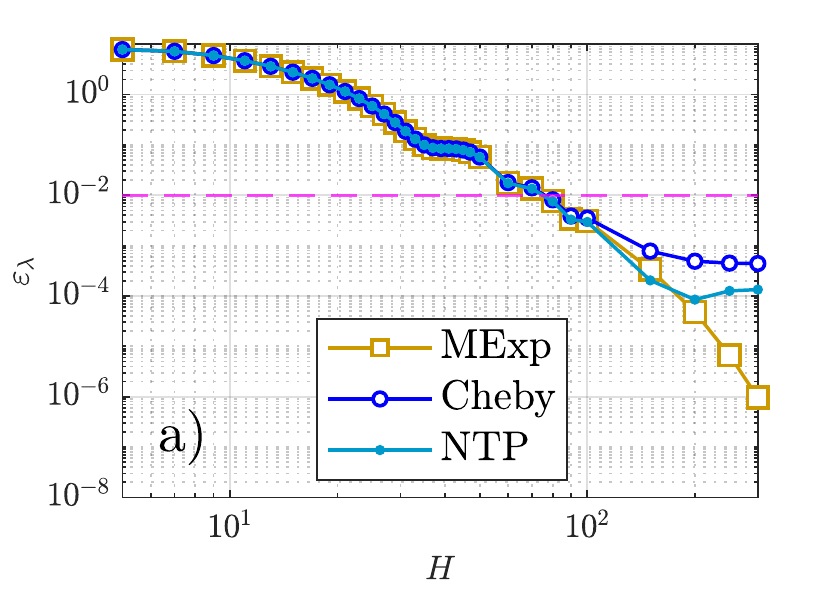}
%\subcaption{}
\end{subfigure}
\hspace{1mm}
\begin{subfigure}[c]{0.485\textwidth}
\includegraphics[width=1\textwidth]{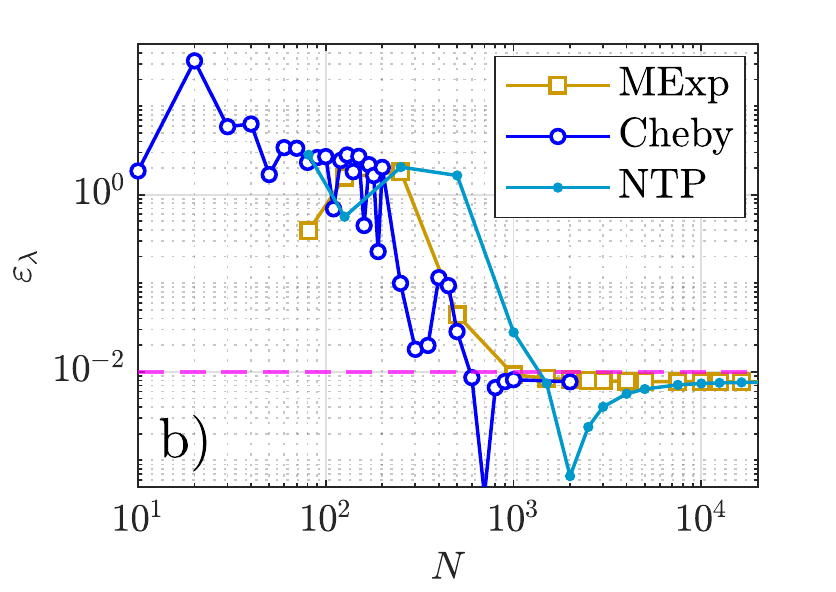}
%\subcaption{}
\end{subfigure}
\caption{Two-degree-of-freedom oscillator with elastic stop: Leading Floquet multiplier error (a) vs. $H$, (b) vs. $N$ for the point with $\Omega/\omega_1=1.2$ indicated in \fref{2DOFunilateral}
}
\label{fig:2DOFunilateral_Hconv}
\end{figure*}
%
% $N=2^{14}$ sufficiently large to avoid visible effect on depicted results % time resolution analogous

%%%%%%%%%%%%%%%%%%%%%%%%%%%%%%%%%%%%%%%%%%%%%%%%%%%%%%%%%%%%%%%%%%%%%%%%%%%%%%%%%%%%%%%%%%%%%%%
\section{Numerical results: Error bound\label{sec:numericalResultsError}}
Two numerical examples are considered in this section:
First, the forced-damped Duffing oscillator with softening characteristic is analyzed, where we shall rigorously prove the existence of some isolated solutions.
Second, the two-degree-of-freedom oscillator with elastic stop from \ssref{2DOFunilateral_benchmark} is revisited.

\subsection{Duffing oscillator with softening characteristic\label{sec:errorestimation_duffing}}
The Duffing oscillator is the archetype of a nonlinear dynamical system.
In the case of a negative coefficient of the cubic spring (softening characteristic), an isolated branch of frequency responses, emerging from zero excitation frequency, is commonly depicted in textbooks.
In \cite{Wagner.2019}, it has been questioned whether such isolated solutions actually exist.
We consider the specific equation of motion,
\ea{
\ddot q + 0.12 \dot q + q - 0.1 q^3 = 0.2 \cos(\Omega t)\fp
 \label{eq:Duffing}
}
\fref{Duffing} shows the amplitude-frequency curve.
As expected for the softening characteristic, the main branch is bent towards the left giving rise to two turning points.
Moreover, an isolated branch occurs.
In fact, multiple of such isolated branches appear as solutions of the HB equations \cite{Wagner.2019}; in this work, we study only a representative one.
While there is a scientific debate concerning the existence of an exact periodic solution, there is consensus that if the solution exists, it will be unstable.
Hence, the existence of these solutions is of rather academic interest.
\begin{figure*}[h]
\centering
\begin{subfigure}[c]{0.485\textwidth}
\includegraphics[width=1\textwidth]{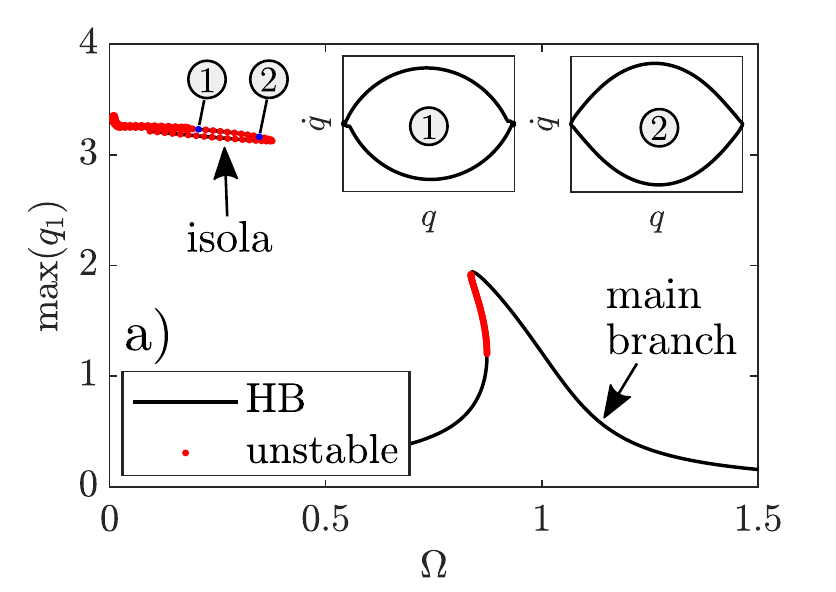}
%\subcaption{}
\end{subfigure}
\hspace{1mm}
\begin{subfigure}[c]{0.485\textwidth}
\includegraphics[width=1\textwidth]{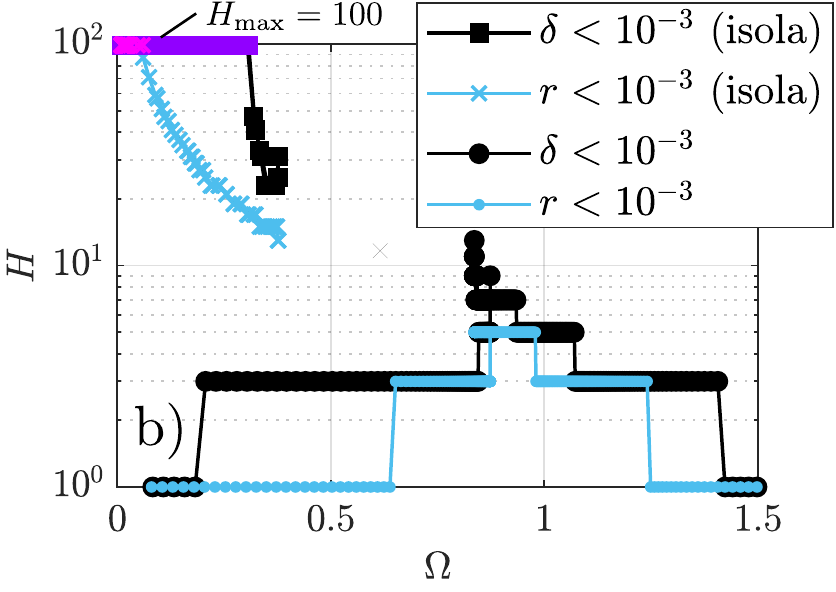}
%\subcaption{}
\end{subfigure}
\caption{Duffing oscillator with softening characteristic: (a) amplitude-frequency curve, (b) harmonic order obtained by adaptive refinement during path.
 {
If the adaptive $H$-refinement reaches the maximum value ($H=100$), the existence of an exact solution could not be proven within the error threshold (indicated by 'pink' and 'purple' markers).
}
}
% (a) $H=9$, $\naft = 2^{13}$, $C=200$
% (b) $C$ adapted with $H$
\label{fig:Duffing}
\end{figure*}
\\
% METHODOLOGICAL DETAILS
The method described in \sref{errorestimation_method} is now applied in order to analyze the existence of an exact periodic solution and to obtain an error bound.
To this end, an upper bound, $\Delta(\delta)$, for the variation of the Jacobian in the $\delta$-neighborhood of the $H$-order approximation, $q_H$, is needed, so that
\ea{
\Delta(\delta) {\leq} \max_{\tau} \left| -0.1 \cdot 3 \cdot \left( \left(q_H\left(\tau\right)+\delta\right)^2 - q_H^2\left(\tau\right)\right)\right|\fp
}
One can see that the expression on the right increases monotonically with $q_H$.
Thus, an upper bound for $q_H(\tau)$ (truncated Fourier series) must be found, which is straight-forward to obtain.
The number of samples within the Alternating Frequency-Time scheme is set to $\naft=4H+1$ to ensure that aliasing is avoided.
The time resolution for Chebyshev interpolation is set to $C=5H$, which was found sufficient to ensure that the depicted results do not visibly change for a further increase of $C$.
\\
% ADAPTIVE H-REFINEMENT; EXISTENCE RESULTS
The computation of the error bound is used in an adaptive $H$-refinement, embedded into a numerical path continuation scheme.
The algorithm reads as follows: % To determine the minimum $H$ required to meet the accuracy requirement,
%
%\begin{enumerate}
%  \item Start the continuation with $H=1$.
%  \item Compute $r$, $M$ as proposed in \sref{errorestimation_method}.
%  \item If a finite error bound $\delta$ can be given and it holds that $\delta<10^{-3}$, go to 5. Otherwise go to 4.
%  \item Increase $H\rightarrow H+2$. If $H=H_{\max}$, go to 5. Otherwise go to 2.
%  \item Do the next step of the continuation.
%\end{enumerate}
%
% An upper limit of $H_{\max}=100$ is specified.
%
\begin{enumerate}
  \item Start the continuation with $H=1$.
  \item Compute $r$, $M$ as proposed in \sref{errorestimation_method}.
  \item If a finite error bound $\delta$ can be given, and $H$ was increased at least once such that $\delta(H-2)>10^{-3}$ and $\delta(H)\leq 10^{-3}$, go to 5. Otherwise go to 4.
  \item If $\delta>10^{-3}$: If $H=H_{\max}$, go to 5. Otherwise increase $H\rightarrow H+2$ and go to 2.\\
  Otherwise ($\delta \leq 10^{-3}$): If $H=H_{\min}$, go to 5. Otherwise decrease $H\rightarrow H-2$ and go to 2.
  \item Do the next step of the continuation.
\end{enumerate}
The limits of $H_{\max}=100$ and $H_{\min}=1$ are specified.
 {
Using this algorithm, we obtain the minimum harmonic truncation order, $H$, required to satisfy a certain convergence criterion with respect to $\delta$.
}
The results are shown in \fref{Duffing}b, both for the main and the isolated branch.
Additionally, results are shown for the case where the convergence criterion in step 3-4 of the above algorithm is replaced by the simple condition $r<10^{-3}$.
It should be emphasized that $r$ and $\delta$ have quite different physical meanings: While $r$ is the (time-domain) residual of a dynamic force balance, $\delta$ is an error bound between approximation and exact solution in the state space.
Thus, using the same numerical threshold $10^{-3}$ is arbitrary and has no justification.
 {
Regarding the number of required $H$, see \fref{Duffing}b,
} 
the condition $r<10^{-3}$ is weaker than $\delta<10^{-3}$; \ie, $r<10^{-3}$ does not imply $\delta<10^{-3}$, but it can be confirmed that  {HB approximations} that satisfy $\delta<10^{-3}$ all have $r<10^{-3}$.
%As can be seen from \fref{Duffing}b, the condition $r<10^{-3}$ is weaker than $\delta<10^{-3}$ as smaller $H$ required to satisfy this criterion; \ie, $r<10^{-3}$ does not imply $\delta<10^{-3}$, but it can be confirmed that the points that satisfy $\delta<10^{-3}$ all have $r<10^{-3}$.
 {
It should further be noted, that the actual error bound $\delta$ (or time domain residual $r$) may be well below the specified threshold value.
}
For some points near the tip of the isolated branch, the existence of a periodic solution in the corresponding $\delta$-neighborhood can be rigorously proven {, associated with an harmonic truncation order $H_\mathrm{min}\leq H < H_\mathrm{max}$}.
As $\Omega\to0$, both $H$-refinement algorithms run into the limit $H=H_{\max}$, leading to inconclusive existence results in a certain range  {(indicated by 'pink' and 'purple' markers)}.
\begin{figure*}[h]
\centering
\begin{subfigure}[c]{0.485\textwidth}
\includegraphics[width=1\textwidth]{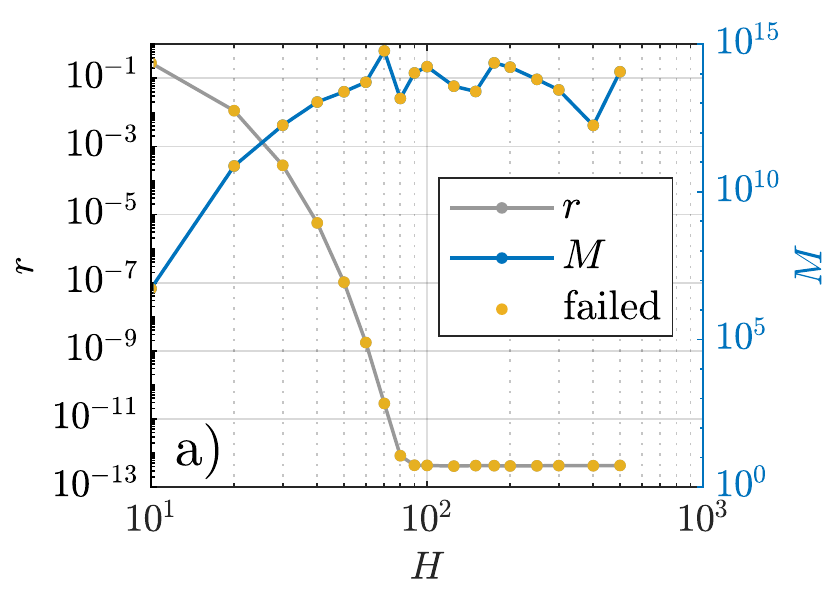}
%\subcaption{}
\end{subfigure}
\hspace{1mm}
\begin{subfigure}[c]{0.485\textwidth}
\includegraphics[width=1\textwidth]{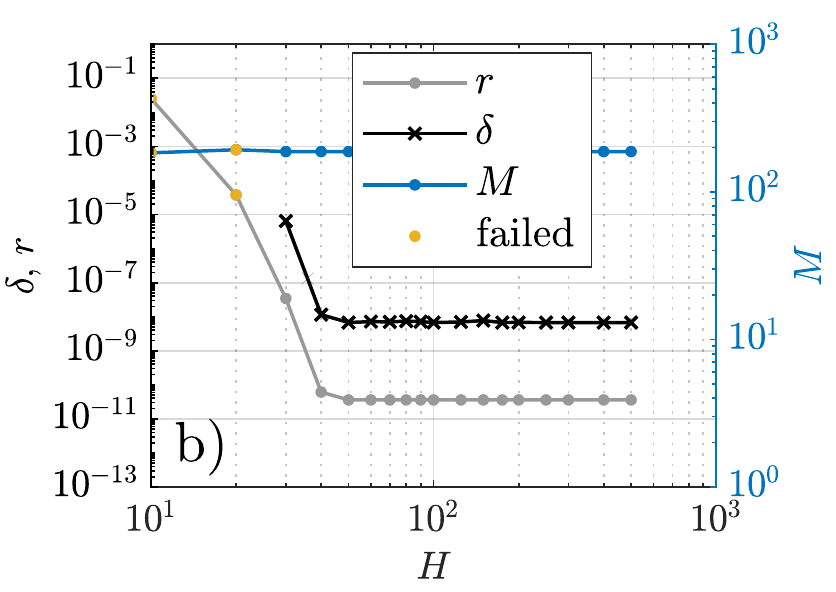}
%\subcaption{}
\end{subfigure}
\caption{Duffing oscillator with softening characteristic: $H$-convergence for two points on the isolated branch (a) $\Omega=0.2$, (b) $\Omega = 0.35$.
 {
For the points indicated as failed ('orange' markers), no error bound could be determined for the given harmonic order.
}
}
\label{fig:Duffing_convergence}
\end{figure*}
% $N=2^{13}$ and $N_\mathrm{Cheby}=C=5\cdot H$
\\
% H-CONVERGENCE
\fref{Duffing_convergence} shows the $H$-convergence for the two different points on the isolated branch indicated in \fref{Duffing}a, for which also the phase projections are shown.
For $\Omega=0.35$, an error bound can be given from $H\geq30$; the error bound is $\delta<10^{-5}$.
The round-off plateau of $r$ is reached at about $H=50$.
Interestingly, the parameter $M$, which is a measure for the error propagation along the cycle, does not change significantly with $H$.
Thus, Urabe's error bound is here mainly driven by the time-domain residual $r$.
For $\Omega=0.2$, $M$ is extremely large ($\mathcal O(M) \approx 10^{14}$ vs. $\mathcal O(M) \approx 10^{2}$ for $\Omega = 0.35$), so that no error bound can be given, even when the round-off plateau of $r$ is reached.
The high $M$ is attributed to the strongly repelling nature of the corresponding region in the state space.
Probably, a higher number precision is needed to obtain an error bound and rigorously prove the existence of an exact solution in the associated $\delta$-neighborhood.
The phase projections in \fref{Duffing}a suggest that the cycle is close to a hetero-clinic orbit, for which HB is known to show poor convergence.

\subsection{Two-degree-of-freedom oscillator with elastic stop \label{sec:errorestimation_unilateral}}
We now revisit the example of \ssref{2DOFunilateral_benchmark}.
Analogously to the previous example, we have as upper bound for the variation of the Jacobian,
\ea{
\Delta(\delta){\leq} \max_{\tau} \left| ~\left.\frac{\dd f_{\mathrm{nl}}}{\dd q_1}\right|_{q_{1,H}(\tau)+\delta} ~-~ \left.\frac{\dd f_{\mathrm{nl}}}{\dd q_1}\right|_{q_{1,H}(\tau)}~\right|\fk\\
\text{where}\quad \frac{\dd f_{\mathrm{nl}}}{\dd q_1} = \frac{100}{2} \left(1 + \frac{100(q_1-1)}{\sqrt{100^2(q_1-1)^2 + 4\varepsilon_\mathrm{reg} }}\right)\fp
}
Continuation in conjunction with adaptive $H$-refinement is carried out as in the previous example, only the limit was raised to $H_{\max} = 500$.
A number of $\naft=2^{13}$ samples per period were used within the Alternating Frequency-Time scheme, and an adaptive time resolution with $C=8H$ was used for Chebyshev interpolation.
\fref{NDOFunilateral_H_vs_omega}a shows results obtained during continuation of the branch depicted in \fref{2DOFunilateral}a.
As in the previous example, higher orders $H$ are needed to obtain an error bound $\delta<10^{-3}$ than to reduce the residual below $r<10^{-3}$.
In fact, $H=500$ is insufficient near the resonance peak  {(indicated by 'pink' markers)}.
For the point with $\Omega/\omega_1=1.2$ indicated in \fref{2DOFunilateral}a, $H>1000$ is needed to obtain an error bound (\fref{NDOFunilateral_H_vs_omega}b).
\begin{figure*}[h]
\centering
\begin{subfigure}[c]{0.475\textwidth}
\includegraphics[width=1\textwidth]{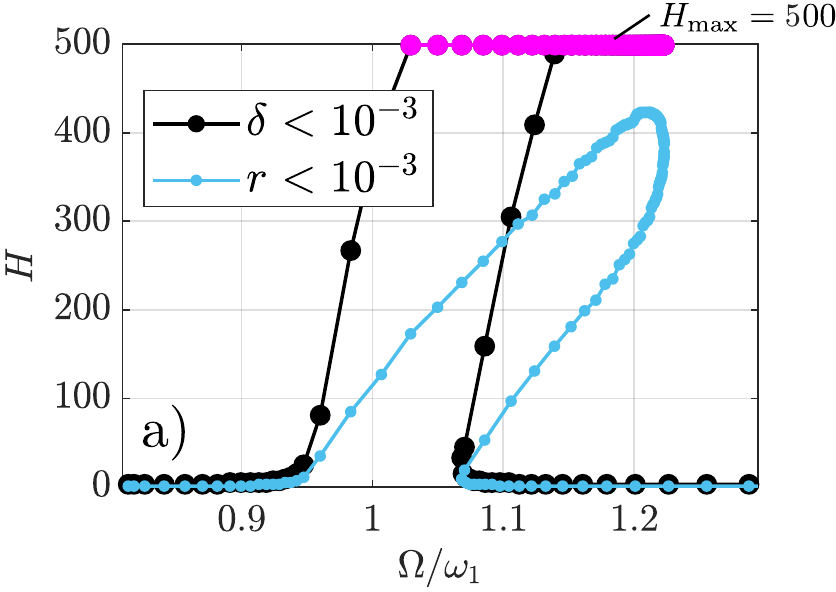}
%\subcaption{}
\end{subfigure}
\hspace{3mm}
\begin{subfigure}[c]{0.475\textwidth}
\includegraphics[width=1\textwidth]{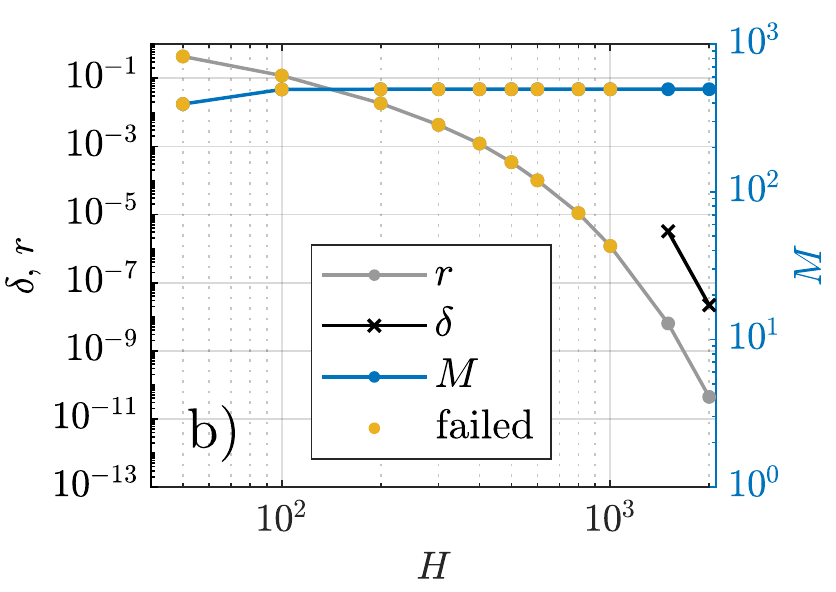}
%\subcaption{}
\end{subfigure}
\caption{Two-degree-of-freedom oscillator with elastic stop: (a) harmonic order obtained by adaptive $H$-refinement during path continuation.
 {
If it reaches the maximum value ($H=500$), the existence of an exact solution could not be proven within the error threshold (indicated by 'pink' markers).
} 
(b) $H$-convergence for the point with $\Omega/\omega_1=1.2$ indicated in \fref{2DOFunilateral}a.
 {
For the points indicated as failed ('orange' markers), no error bound could be determined for the given harmonic order.
} 
}
\label{fig:NDOFunilateral_H_vs_omega}
\end{figure*}

\section{Conclusions \label{sec:conclusions}}
% OVERVIEW
The purpose of the present work was to determine whether (a) Chebyshev-based stability analysis and (b) Urabe's error bound are useful features for Harmonic Balance.
% CHEBYSHEV
The idea of the Chebyshev-based stability analysis is to exploit the properties of Chebyshev polynomials to avoid the time step integration of the monodromy matrix.
Indeed, the time integration is replaced by the solution of a single linear algebraic equation system of dimension $dN$ where $d$ is the number of degrees of freedom and $N$ is the number of time levels per period (equal to the number of Chebyshev polynomials).
Thus, the method scales quadratically with $N$, in contrast to conventional time step integration (\eg Newmark or Matrix Exponential Method), which scales linearly with $N$.
Consequently, the Chebyshev Method tends to be quicker when only a relatively small $N$ is needed.
This is usually the case when also a low harmonic truncation order $H$ is sufficient, as is typically the case for (inherently) smooth nonlinearities.
In contrast, in the presence of (sharply regularized) discontinuities, high $H$ and $N$ are typically needed, and the method becomes inefficient.
Reasoning that for high $H$, HB itself becomes inefficient compared to alternative methods for computing periodic oscillations, the Chebyshev-based stability analysis is regarded as natural supplement.
About one  {to two orders} of magnitude speedup compared to the state-of-the-art alternatives has been found in the corresponding numerical examples.
It should be stressed that further potential for computational speedup exists.
For the case of local nonlinearities, in particular, the computational burden of solving the large-dimensional equation system could be reduced by a block condensation with pre-computed factorization of the linear (time-invariant) part.
\\
% URABE
Concerning Urabe's error bound, a similar distinction should be made with regard to the considered problem class.
For inherently smooth problems, Urabe's theorem appears useful in computational engineering practice.
In contrast, in the presence of (sharply regularized) discontinuities, an extremely high truncation order $H$ is needed to reduce the residual sufficiently to obtain a finite error bound.
In that case, conventional time step integration of the periodic solution may become faster so that the use of HB becomes somewhat obsolete.
As soon as an error bound could be given, for the considered numerical examples, this error bound was so tight that a further increase of the truncation order never seemed reasonable.
From a practical perspective, it would be desirable to have a rigorous existence proof and a finite error bound already for much lower truncation orders.
Further, the computational effort for obtaining the measures needed for Urabe's theorem is relatively high, even when Chebyshev polynomials are again used to achieve numerical efficiency.
More specifically, the bottleneck is the computation of the measure $M$, which gives an upper bound for the error propagation within the period.
Given the popularity of HB in engineering, research towards the mathematically rigorous development of computationally cheaper error bounds is desirable.

\section*{Acknowledgements}
This work was funded by the Deutsche Forschungsgemeinschaft (DFG, German Research Foundation) [Project 438529800].

\appendix
\section{Chebyshev analysis \label{asec:Cheby}}
%\COMMENT{PLEASE CHECK. ON FIRST SIGHT, IT LOOKS LIKE A CONTRADICTION TO WIKI.}
%\\
% INTEGRATION RULE
According to the \emph{integration rule} of Chebyshev polynomials,
\ea{
b_C &= \int a_C\dd\tau \fk \label{eq:integrationruleone}\\
\check{\mm b} &= \mm G \check{\mm a}\fk \label{eq:integrationruletwo}
}
\ie; the integral of a Chebyshev polynomial is a Chebyshev polynomial of the same order.
The coefficients are obtained as expressed in \eref{integrationruletwo}, using an integration operational matrix $\mm G$, which reads
\ea{
\mms G=\pi
\begin{bmatrix}
1 & 1 & 0 & 0   & \hdots & 0 \\
-1/4 & 0 & 1/4 & 0 & & 0 \\
-1/3 & -1/2 & 0 & 1/6  & & 0\\
1/8 & 0 & -1/4 & 0 &  & 0\\
\vdots & & & &  \ddots  & \dfrac{1}{2(C-1)}\\
 \dfrac{-(-1)^{C-1}}{C(C-2)} & & & &  \dfrac{-1}{2(C-2)}  & 0
\end{bmatrix},
%\begin{bmatrix}
%1/2     & -1/8      & -1/6  & 1/16  & \hdots                        & \dfrac{-(-1)^{C-1}}{2C(C-2)} \\
%1/2     & 0             & -1/4 & 0          &                                   &  \\
%0         & 1/8        & 0        & -1/8    &                                   & \\
%0         & 0           & 1/12  & 0             &                                   &\\
%0         & 0            & 0        & 1/16      &                                   &  \\
%\vdots &            &              &                &     \ddots              &  \dfrac{-1}{4(C-2)} \\
%0       & 0             & 0         & 0             & \dfrac{1}{4(C-1)} & 0
%\end{bmatrix}\fk
}
for the considered case of Chebyshev base functions of the first kind.
The factor $\pi$ stems from the considered interval length $[0,2\pi]$ rather than the unit interval $[-1,1]$.
%where $\square^{\mathrm T}$ denotes transpose.
\\
% PRODUCT RULE
According to the \emph{product rule} of Chebyshev polynomials, the product, $c=a_C b_C$ of two Chebyshev polynomials $a_C$ and $b_C$ is a Chebyshev polynomial of higher order.
When only the coefficients up to order $C$ are of interest, these can be expressed as
\ea{
\check{\mm c} = \mm P\left(\check{\mm a}\right) \check{\mm b}\fk \label{eq:productrule}
}
using the product operational matrix
\ea{
\mm P\left(\check{\mm a}\right) =
\begin{bmatrix}
\check a(1) & \check a(2)/2 & \check a(3)/2  & \hdots & \check a(C)/2 \\
\check a(2) & \check a(1) + \check a(3)/2 & ( \check a(2) +\check a(4) )/2 & & \check a(C-1)/2 \\
\check a(3) &  ( \check a(2) +\check a(4) )/2 & \check a(1) + \check a(5)/2 &  & \check a(C-2)/2 \\
\vdots & & & \ddots  & \\
 \check a(C) & \check a(C-1)/2  &  &  & \check a(1)
\end{bmatrix}\fp
}
\section*{References}
\bibliographystyle{elsarticle-num}
\bibliography{literature}

\begin{thebibliography}{10}
\expandafter\ifx\csname url\endcsname\relax
  \def\url#1{\texttt{#1}}\fi
\expandafter\ifx\csname urlprefix\endcsname\relax\def\urlprefix{URL }\fi
\expandafter\ifx\csname href\endcsname\relax
  \def\href#1#2{#2} \def\path#1{#1}\fi

\bibitem{urab1965}
M.~Urabe, Galerkin's procedure for nonlinear periodic systems, Archive for
  Rational Mechanics and Analysis 20~(2) (1965) 120--152.

\bibitem{came1989}
T.~M. Cameron, J.~H. Griffin, An alternating frequency/time domain method for
  calculating the steady-state response of nonlinear dynamic systems, Journal
  of Applied Mechanics 56~(1) (1989) 149--154.
\newblock \href {https://doi.org/10.1115/1.3176036}
  {\path{doi:10.1115/1.3176036}}.

\bibitem{padm1995b}
C.~Padmanabhan, R.~Singh, Analysis of periodically excited non-linear systems
  by a parametric continuation technique, Journal of Sound and Vibration
  184~(1) (1995) 35--58.
\newblock \href {https://doi.org/10.1006/jsvi.1995.0303}
  {\path{doi:10.1006/jsvi.1995.0303}}.

\bibitem{Krack.2019}
M.~Krack, J.~Gross, Harmonic Balance for Nonlinear Vibration Problems,
  Springer, 2019.
\newblock \href {https://doi.org/10.1007/978-3-030-14023-6}
  {\path{doi:10.1007/978-3-030-14023-6}}.

\bibitem{Hall.2002}
K.~C. Hall, J.~P. Thomas, W.~S. Clark, Computation of unsteady nonlinear flows
  in cascades using a harmonic balance technique, AIAA Journal 40~(5) (2002)
  879--886.

\bibitem{Gilmore.1991}
R.~J. Gilmore, M.~B. Steer, Nonlinear circuit analysis using the method of
  harmonic balance---a review of the art. part i. introductory concepts,
  International Journal of Microwave and Millimeter--Wave Computer--Aided
  Engineering 1~(1) (1991) 22--37.

\bibitem{Hsu.1972}
C.~S. Hsu, Impulsive parametric excitation: theory, Journal of Applied
  Mechanics 39~(2) (1972) 551--558.

\bibitem{Friedmann.1977}
P.~Friedmann, C.~E. Hammond, T.-H. Woo, Efficient numerical treatment of
  periodic systems with application to stability problems, International
  Journal for Numerical Methods in Engineering 11 (1977) 1117--1136.

\bibitem{peletan2013}
L.~Peletan, S.~Baguet, M.~Torkhani, G.~Jacquet-Richardet, A comparison of
  stability computational methods for periodic solution of nonlinear problems
  with application to rotordynamics, Nonlinear Dynamics 72~(3) (2013) 671--682.

\bibitem{groll2001b}
G.~v. Groll, D.~J. Ewins, The harmonic balance method with arc-length
  continuation in rotor/stator contact problems, Journal of Sound and Vibration
  241~(2) (2001) 223--233.
\newblock \href {https://doi.org/10.1006/jsvi.2000.3298}
  {\path{doi:10.1006/jsvi.2000.3298}}.

\bibitem{laza2010}
A.~Lazarus, O.~Thomas, A harmonic-based method for computing the stability of
  periodic solutions of dynamical systems, Comptes Rendus M{\'e}canique 338~(9)
  (2010) 510--517.

\bibitem{Moore.2005}
{Gerald Moore}, Floquet theory as a computational tool, SIAM Journal on
  Numerical Analysis 42~(6) (2005) 2522--2568.
\newblock \href {https://doi.org/10.1137/S0036142903434175}
  {\path{doi:10.1137/S0036142903434175}}.

\bibitem{Wagner.2019}
U.~von Wagner, L.~Lentz, On the detection of artifacts in harmonic balance
  solutions of nonlinear oscillators, Applied Mathematical Modelling 65 (2019)
  408--414.
\newblock \href {https://doi.org/10.1016/j.apm.2018.08.013}
  {\path{doi:10.1016/j.apm.2018.08.013}}.

\bibitem{Wagner.2018}
U.~von Wagner, L.~Lentz, On artifact solutions of semi-analytic methods in
  nonlinear dynamics, Archive of Applied Mechanics 88~(10) (2018) 1713--1724.
\newblock \href {https://doi.org/10.1007/s00419-018-1397-3}
  {\path{doi:10.1007/s00419-018-1397-3}}.

\bibitem{Grolet.2012}
A.~Grolet, F.~Thouverez, On a new harmonic selection technique for harmonic
  balance method, Mechanical Systems and Signal Processing 30 (2012) 43--60.
\newblock \href {https://doi.org/10.1016/j.ymssp.2012.01.024}
  {\path{doi:10.1016/j.ymssp.2012.01.024}}.

\bibitem{jaum2010}
V.~Jaumouill{\'e}, J.~J. Sinou, B.~Petitjean, An adaptive harmonic balance
  method for predicting the nonlinear dynamic responses of mechanical systems -
  application to bolted structures, Journal of Sound and Vibration 329~(19)
  (2010) 4048--4067.

\bibitem{Ferri.2009}
A.~Ferri, M.~Leamy, Error estimates for harmonic-balance solutions of nonlinear
  dynamical systems, in: 50th AIAA/ASME/ASCE/AHS/ASC Structures, Structural
  Dynamics, and Materials Conference, Structures, Structural Dynamics, and
  Materials and Co-located Conferences, {American Institute of Aeronautics and
  Astronautics}, 2009.
\newblock \href {https://doi.org/10.2514/6.2009-2667}
  {\path{doi:10.2514/6.2009-2667}}.

\bibitem{Stokes.1972}
A.~Stokes, On the approximation of nonlinear oscillations, Journal of
  Differential Equations 12~(3) (1972) 535--558.
\newblock \href {https://doi.org/10.1016/0022-0396(72)90024-1}
  {\path{doi:10.1016/0022-0396(72)90024-1}}.

\bibitem{urab1966}
M.~Urabe, A.~Reiter, Numerical computation of nonlinear forced oscillations by
  galerkin's procedure, Journal of Mathematical Analysis and Applications
  14~(1) (1966) 107--140.

\bibitem{vanDooren.1988}
R.~{van Dooren}, On the transition from regular to chaotic behaviour in the
  duffing oscillator, Journal of Sound and Vibration 123~(2) (1988) 327--339.
\newblock \href {https://doi.org/10.1016/S0022-460X(88)80115-9}
  {\path{doi:10.1016/S0022-460X(88)80115-9}}.

\bibitem{Baker.2005}
A.~Baker, M.~Dellnitz, O.~Junge, Topological method for rigorously computing
  periodic orbits using fourier modes, Discrete and Continuous Dynamical
  Systems - DISCRETE CONTIN DYN SYST 13 (2005).
\newblock \href {https://doi.org/10.3934/dcds.2005.13.901}
  {\path{doi:10.3934/dcds.2005.13.901}}.

\bibitem{Lessard.2014}
J.-P. Lessard, C.~Reinhardt, Rigorous numerics for nonlinear differential
  equations using chebyshev series, SIAM J. Numer. Anal. 52 (2014) 1--22.
\newblock \href {https://doi.org/10.1137/13090883X}
  {\path{doi:10.1137/13090883X}}.

\bibitem{Lessard.2015}
J.-P. Lessard, J.~{Mireles James}, A.~Hungria, Rigorous numerics for analytic
  solutions of differential equations: The radii polynomial approach,
  Mathematics of Computation 85 (2015).
\newblock \href {https://doi.org/10.1090/mcom/3046}
  {\path{doi:10.1090/mcom/3046}}.

\bibitem{Kogelbauer.2021}
F.~Kogelbauer, T.~Breunung, When does the method of harmonic balance give a
  correct prediction for mechanical systems?, Applicable Analysis 2 (2021)
  1--19.
\newblock \href {https://doi.org/10.1080/00036811.2021.1953482}
  {\path{doi:10.1080/00036811.2021.1953482}}.

\bibitem{boydspectral}
J.~P. Boyd, Chebyshev and Fourier spectral methods, {Courier Corporation},
  2001.

\bibitem{sinh1991}
S.~C. Sinha, D.~H. Wu, An efficient computational scheme for the analysis of
  periodic systems, Journal of Sound and Vibration 151~(1) (1991) 91--117.
\newblock \href {https://doi.org/10.1016/0022-460X(91)90654-3}
  {\path{doi:10.1016/0022-460X(91)90654-3}}.

\bibitem{Wu.1994}
D.-H. Wu, S.~C. Sinha, A new approach in the analysis of linear systems with
  periodic coefficients for applications in rotorcraft dynamics, The
  Aeronautical Journal (1968) 98~(971) (1994) 9--16.
\newblock \href {https://doi.org/10.1017/S0001924000050302}
  {\path{doi:10.1017/S0001924000050302}}.

\bibitem{Woiwode.2020}
L.~Woiwode, N.~N. Balaji, J.~Kappauf, F.~Tubita, L.~Guillot, C.~Vergez,
  B.~Cochelin, A.~Grolet, M.~Krack, Comparison of two algorithms for harmonic
  balance and path continuation, Mechanical Systems and Signal Processing 136
  (2020) 106503.
\newblock \href {https://doi.org/10.1016/j.ymssp.2019.106503}
  {\path{doi:10.1016/j.ymssp.2019.106503}}.

\bibitem{GarciaSaldana.2013}
J.~D. Garc{\'i}a-Salda{\~n}a, A.~Gasull, A theoretical basis for the harmonic
  balance method, Journal of Differential Equations 254~(1) (2013) 67--80.
\newblock \href {https://doi.org/10.1016/j.jde.2012.09.011}
  {\path{doi:10.1016/j.jde.2012.09.011}}.

\bibitem{Thouverez.2003}
F.~Thouverez, Presentation of the ecl benchmark, Mechanical Systems and Signal
  Processing 17~(1) (2003) 195--202.
\newblock \href {https://doi.org/10.1006/mssp.2002.1560}
  {\path{doi:10.1006/mssp.2002.1560}}.

\bibitem{Noel.2015}
J.-P. Noël, L.~Renson, C.~Grappasonni, G.~Kerschen, Identification of
  nonlinear normal modes of engineering structures under broadband forcing,
  Mechanical Systems and Signal Processing 74 (2016) 95--110, special Issue in
  Honor of Professor Simon Braun.
\newblock \href {https://doi.org/https://doi.org/10.1016/j.ymssp.2015.04.016}
  {\path{doi:https://doi.org/10.1016/j.ymssp.2015.04.016}}.

\end{thebibliography}

\end{document}